\documentclass[11pt]{amsart}
\usepackage{amsmath,amsthm, epsf, amscd}
\usepackage{times}

\usepackage[dvips]{graphicx}
\usepackage{psfrag}
\usepackage[psamsfonts]{amssymb}
\usepackage[all]{xy}

\numberwithin{equation}{section}

\newtheorem{Thm}{Theorem}[section]
\newtheorem{Prop}[Thm]{Proposition}
\newtheorem{Lem}[Thm]{Lemma}
\newtheorem{Cor}[Thm]{Corollary}

\newtheorem{Obs}[Thm]{Observation}

\theoremstyle{remark}
\newtheorem{Rem}[Thm]{Remark}
\newtheorem{Exa}[Thm]{Example}

\theoremstyle{definition}
\newtheorem{Def}[Thm]{Definition}
\newtheorem{Prob}[Thm]{Problem}

\newcommand{\fig}[1]
        {\raisebox{-0.5\height}
                 {\includegraphics{#1}}
        }

\def\Z{{\mathbb Z}}
\def\R{{\mathbb R}}
\def\Q{{\mathbb Q}}
\def\C{{\mathbb C}}
\def\calA{\mathcal{A}}
\def\calD{\mathcal{D}}

\def\calG{\mathcal{G}}

\def\calV{\mathcal{V}}

\def\Im{\mathrm{Im}\,}
\def\deg{\mathrm{deg}\,}
\def\span{\mathrm{span}}
\def\Hom{\mathrm{Hom}}
\def\eqdef{\stackrel{\rm def}{=}}
\def\Aut{\mathrm{Aut}\,}
\def\inn{\mathrm{int}\,}
\def\cM{\overline{M}}

\def\BDiff{\widetilde{B\mathrm{Diff}}}
\def\EDiff{\widetilde{E\mathrm{Diff}}}
\def\Emb{\mathrm{Emb}}
\def\Diff{\mathrm{Diff}}
\def\gtimes{\ltimes}

\newcommand{\mapright}[1]{
	\smash{\mathop{
		\hbox to 1cm{\rightarrowfill}}\limits^{#1}}}
\newcommand{\mapsright}[1]{
	\stackrel{#1}{\to}}
\newcommand{\mapleft}[1]{
	\smash{\mathop{
		\hbox to 1cm{\leftarrowfill}}\limits^{#1}}}

\begin{document}
\title[On Kontsevich's characteristic classes]{On Kontsevich's characteristic classes for smooth 5- and 7-dimensional homology sphere bundles}

\author[T. Watanabe]{Tadayuki Watanabe}
\address{Research Institute for Mathematical Sciences, Kyoto University, Kyoto, Japan}
\email{tadayuki@kurims.kyoto-u.ac.jp}

\date{\today}

\begin{abstract}
M. Kontsevich constructed universal characteristic classes of smooth bundles with fiber a framed odd-dimensional integral homology sphere. In dimension 3, they are known to give a universal finite type invariants of homology 3-spheres. However, they have not been well understood for higher fiber dimensions. The purpose of the present paper is twofold. First, we obtain a bordism invariant of smooth unframed bundles with fiber a 5-dimensional homology sphere, which is defined as a sum of the simplest Kontsevich class and the second signature defect. It may be in some sense a higher dimensional analogue of the Casson invariant. Second, we construct a family of $M$-bundles. By evaluating on those $M$-bundles, we show that Kontsevich's universal characteristic classes are highly non-trivial in the case of fiber dimension 7. As a corollary, new estimates for unstable rational homotopy groups of $\mathrm{Diff}(D^7\mbox{ rel }\partial)$ are obtained.
\end{abstract}
\maketitle
%\tableofcontents
%%%%%%%%%%%%%%%%%%%%%%%%%%%%%%
%%%%%%%%%%%%%%%%%%%%%%%%%%%%%%
%%%%%%%%%%%%%%%%%%%%%%%%%%%%%%
\section{Introduction}

In \cite{Kon}, M. Kontsevich constructed universal characteristic classes of smooth framed $M$-bundles with fiber an odd dimensional integral homology sphere $M$. The construction of the Kontsevich classes involves the graph complex and configuration space integrals (see \cite{Kon} or \S\ref{ss:kon-cc} for the definition). In the case of 3-dimensional homology spheres, the Kontsevich classes are 0-forms, i.e., real valued diffeomorphism invariants, and it is shown in \cite{KT} that all the Ohtsuki finite type invariants (\cite{Oh}) are recovered in this way. It is also known that there are very many Ohtsuki finite type invariants, hence the Kontsevich classes for 3-dimensional homology spheres are very strong.

In the present paper, we study the Kontsevich classes for higher odd-dimensional homology spheres. In particular, we get some results for 5- and 7-dimensional homology spheres. There are roughly two parts in the present paper, each of which can be read almost separately.

First, the Kontsevich classes are the characteristic classes for smooth `framed' $M$-bundles. Let $M^\bullet$ denote $M$ with a puncture at a fixed point $\infty\in M$. By a framing on an $M$-bundle, we mean a trivialization of $TM^\bullet$ along each fiber that is standard near $\partial M^\bullet$, namely it looks near $\partial M^\bullet$ like the standard Euclidean plane which may be identified as $(S^m)^\bullet$. In the case of 5-dimensional homology spheres, we study in \S\ref{ss:bordism-inv} the framing dependence of the simplest Kontsevich class associated to the $\Theta$-graph, which is a 2-form on the base space, and we obtain a bordism invariant of unframed $M$-bundles by adding a certain multiple of the second signature defect invariant of Hirzebruch (Theorem~\ref{thm:A}). As a corollary, it follows that the simplest Kontsevich class, as a characteristic class for framed bundles, is non-trivial. Our formula for the bordism invariant is in some sense a higher dimensional analogue of Morita's formula for the Casson invariant in \cite{Mo}. 

Second, in the case of 7-dimensional homology spheres, we construct in \S\ref{ss:clasper-bundles} a family of framed $M$-bundles, which we call graph clasper-bundles, by using higher dimensional claspers. Higher dimensional claspers are introduced in \cite{W} as higher dimensional generalizations of Habiro's claspers in 3-dimension \cite{Hab}. Instead of the Borromean rings in Habiro's graph claspers, we use a `suspension' of the higher dimensional Borromean rings. Then we show that our construction is in some sense dual to the Kontsevich classes (Theorem~\ref{thm:B}). Proof of Theorem~\ref{thm:B} is inspired by Kuperberg-Thurston's proof of the universality of their version of Kontsevich's perturbative invariant \cite{KT} and somewhat philosophically by Cattaneo--Cotta-Ramusino--Longoni\cite{CCL}. As a consequence of Theorem~\ref{thm:B}, it turns out that the Kontsevich classes are highly non-trivial and that there are as many Kontsevich classes as the Ohtsuki finite type invariants, for any fixed 7-dimensional homology spheres. This already implies that there are a lot of smooth framed bundles, and that clasper-bundle surgery can be used effectively to produce a lot of bundles similarly as the case of 3-dimensional homology sphere, while usual surgery along framed links in a higher dimensional manifold seems not so effective unless the manifold is nilpotent \cite{W}. Also, our construction gives some linearly independent elements of the homotopy groups $\pi_{4n}{B\Diff}(D^7\mbox{ rel }\partial)\otimes\Q$ of the base of the universal $(D^7\mbox{ rel }\partial)$-bundle, thus we obtain new non-trivial estimates for the rational homotopy groups $\pi_{4n-1}\Diff(D^7\mbox{ rel }\partial)\otimes\Q$ of the infinite dimensional Lie group $\Diff(D^7\mbox{ rel }\partial)$, the group of diffeomorphisms. This is in some sense unstable informations. By the way, in the stable range, the isomorphism $\pi_{4n}B\mathrm{\Diff}(D^{2m-1}\mbox{ rel }\partial)\otimes\Q\cong\Q\ (2m-1>>4n)$ has been given by Farrell-Hsiang \cite{FH}. 

In \S\ref{s:directions}, we will remark some future directions. We think that the study of cohomology classes of the space of certain link embeddings is a higher dimensional generalization of the study of link invariants in a 3-manifold. Similarly, we think that the study of universal characteristic classes is a higher dimensional generalization of the study of invariants of 3-manifolds. We expect that there is a rich theory for smooth bundles as in the theory of Ohtsuki's finite type invariants of homology 3-spheres and that clasper-bundle surgery gives an important correspondence between the two.

%\clearpage
%%%%%%%%%%%%%%%%%%%%%%%%%%%%%%
%%%%%%%%%%%%%%%%%%%%%%%%%%%%%%
%%%%%%%%%%%%%%%%%%%%%%%%%%%%%%
\section{Kontsevich's universal characteristic classes}\label{ss:kon-cc}

Here we briefly review the definition of Kontsevich's universal characteristic classes.

%%%%%%%%%%%%%%%%%%%%%%%%%%%%%%
%%%%%%%%%%%%%%%%%%%%%%%%%%%%%%
\subsection{Feynman diagrams}
First we define the space $\calA_{2n}$ of trivalent graphs. An {\it orientation} on a trivalent graph $\Gamma$ is a choice of an ordering of three edges incident to each trivalent vertex, considered modulo even number of swappings of the orders. We present the orientation in plane diagrams by assuming that the order of three edges incident to each trivalent vertex is given by anti-clockwise order.

Let $\calG_{2n}$ be the real vector space spanned by all connected trivalent graphs with oriented $2n$ vertices. Let $\calA_{2n}$ be the quotient space of $\calG_{2n}$ by the subspace spanned by the vectors of the following form:
\begin{equation}\label{eq:ihx}
\psfrag{-}[cc][cc]{$-$}
\psfrag{+}[cc][cc]{$+$}
\fig{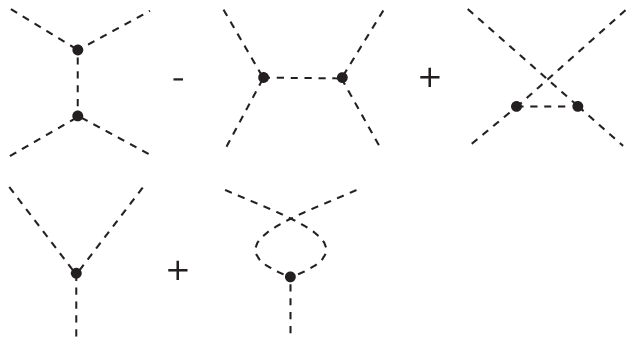}
\end{equation}
We call the vectors in (\ref{eq:ihx}) {\it IHX} and {\it AS relations} respectively. We will write as $[\Gamma]$ the element of $\calA_{2n}$ represented by a graph $\Gamma\in\calG_{2n}$. The {\it degree} of a trivalent graph is defined as the number of vertices. For example, $\calA_2=\span_\R\{[\Theta]\}$, where $[\Theta]$ is the $\Theta$-graph.

%%%%%%%%%%%%%%%%%%%%%%%%%%%%%%
%%%%%%%%%%%%%%%%%%%%%%%%%%%%%%
\subsection{Fulton-MacPherson-Kontsevich compactification of the configuration space}
Let $M$ be an $m$-dimensional homology sphere with a fixed point $\infty\in M$. Let $C_n(M)$ be the Fulton-MacPherson-Kontsevich compactification (\cite{FM}) of the configuration space 
\[ M^{\times n}\setminus\mbox{(diagonals)}. \]
Here we include in the diagonals the set of configurations with some points go infinity. For example, $C_2(M)$ is obtained from $M\times M$ by blowing up first along $(\infty, \infty)$ and then along the disjoint diagonals
\[ \mbox{(the diagonal)}\cup (M\setminus\infty)\times\{\infty\}\cup\{\infty\}\times (M\setminus\infty). \] 
Neighborhood in $C_2(M)$ of the face of $\partial{C}_2(M)$ corresponding to $(\infty,\infty)\in M\times M$ is the same form as a neighborhood of the face of $\partial{C}_2(S^{m-1})$ corresponding to $(\infty,\infty)\in S^{m-1}\times S^{m-1}$, and there exists the Gauss map $p_{S^{m-1}}$ from the $(\infty,\infty)$-face of $\partial{C}_2(S^{m-1})$ to $S^{m-1}$, given by the unit relative vector in $\R^m$. So the Gauss map $p_M$ from the $(\infty,\infty)$-face of $\partial{C}_2(M)$ to $S^{m-1}$ is defined as $p_{S^{m-1}}$. Moreover, the union of the faces corresponding to the above three diagonals is naturally a trivial $S^{m-1}$-bundle. One obtains a map $p_M:\partial{C}_2(M)\setminus((\infty,\infty)\mbox{-face})\to S^{m-1}$ given by the projection onto the $S^{m-1}$-factor. It is determined by the framing. Therefore a continuous map $p_M:\partial{C}_2(M)\to S^{m-1}$ is defined. It is known that $p_M^*\omega_{S^{m-1}}$, where $\omega_{S^{m-1}}$ is the $SO(m)$-invariant unit volume form on $S^{m-1}$,  extends to a closed $(m-1)$-form $\alpha_M$ on $C_2(M)$ and it generates $H^{m-1}(C_2(M);\R)$ \cite{Coh, Les}.

%%%%%%%%%%%%%%%%%%%%%%%%%%%%%%
%%%%%%%%%%%%%%%%%%%%%%%%%%%%%%
\subsection{Universal smooth $M$-bundle}
Let $M^\bullet$ denote $M$ with a puncture at $\infty\in M$. By a {\it smooth vertically framed $M$-bundle}, we mean a smooth bundle with fiber $M$, together with a fixed inclusion $M^\bullet\hookrightarrow M\setminus\{\infty\}$ such that the bundle is trivialized on $M\setminus \mathrm{int}(M^\bullet)$ and such that there is a trivialization of its vertical tangent bundle restricted to $M^\bullet$-fiber, namely, tangent bundle along the $M^\bullet$-fibers, that is also standard near $\partial M^\bullet$. We will call such a trivialization a {\it vertical framing}. 

Let $\widetilde{\Emb}(M\setminus\{\infty\},\R^{\infty})$ be the space of smooth tangentially framed embeddings $M\setminus\{\infty\}\to \R^{\infty}$ that are standard near $\infty$, i.e., coincide with $\R^m\subset\R^{\infty}$ near $\infty$. Here $\R^{\infty}$ denotes the Hilbert space of square summable sequences. We equip $\widetilde{\Emb}(M\setminus\{\infty\},\R^\infty)$ with the $F\calD$-topology in \cite{Mic}. Then the bundle
\[ \pi_{\Diff{M}}:\widetilde{\Emb}(M\setminus\{\infty\},\R^{\infty})\to\widetilde{\Emb}(M\setminus\{\infty\},\R^{\infty})/\mathrm{Diff}(M^\bullet\mbox{ rel }\partial) \]
is a disjoint union of copies of the universal framed $\Diff (M^\bullet\mbox{ rel }\partial)$-bundle\footnote{Here we say universal framed bundle in the sense that it is contractible into the space of framings on $M$ that are standard near $\infty$ and that there is a bijection between the set of isomorphism classes of vertically framed $M$-bundles over $B$ and the homotopy set $[B,\BDiff{M}]$.}, each associated to a homotopy class of framings on $M^\bullet$ (in the case $M^\bullet$ is a punctured homology sphere, there are at most $\Z\times$finite-copies). We denote the bundle $\pi_{\Diff{M}}$ simply by $\EDiff{M}\to\BDiff{M}$. We fix a base point of each component of $\BDiff{M}$ and fix a standard framing on the fiber of it. $\Diff(M^\bullet\mbox{ rel }\partial)$ acts on $\widetilde{\Emb}(M\setminus\{\infty\},\R^{\infty})$ from the right by $((\phi,\tilde\tau_M)\cdot g)(x)=(\phi(gx), \tilde\tau_M(gx))$ for $\phi\in\Emb(M\setminus\{\infty\},\R^\infty)$ and for $\tilde\tau_M:(M^\bullet,\partial M^\bullet)\to (GL_+(\R^m),1)$ being a difference from the standard framing. $\BDiff{M}$ is also considered as the base of the universal smooth framed $M$-bundle 
\[ \pi_M:M\gtimes \EDiff{M}\to\BDiff{M}, \]
 associated to $\pi_{\Diff M}$\footnote{In fact, $\BDiff{M}$ is a kind of an infinite dimensional smooth manifold for which the de Rham theorem holds. See \cite{Mic, Mic2} for details about it.}. Here the expression $F\gtimes \EDiff{M}$ means the Borel construction $F\times_{\Diff(M^\bullet\mbox{\scriptsize \ rel }\partial)}\EDiff{M}$. From general theory of bundles, an isomorphism class of a smooth framed $M$-bundle $E\to B$ is determined by the homotopy class of a classifying map $f:B\to\BDiff{M}$. We will often identify the image of a classifying map $f$ with the induced bundle $f^*\pi_{\Diff M}$ and in the light of this identification we will identify a fiber with a point of $\BDiff{M}$. Usually, cohomology classes of $\BDiff{M}$ are used for homotopy classification of classifying maps and they are called universal characteristic classes (e.g. \cite{Mo2}). For bundles over closed manifolds, bordism invariants $\Omega_*(\BDiff{M})\to\calV$ ($\calV$: a $\Z$-module or a real vector space) may also be used for the classification.
 
From the result of Appendix~\ref{s:alpha-form}, there exists a closed $(m-1)$-form $\alpha_{\Diff{M}}$ on the universal $C_2(M)$-bundle
\[ \pi_{C_2(M)}:C_2(M)\gtimes\EDiff{M}\to \BDiff{M}\]
associated to $\pi_M$, whose restriction on each fiber is $[\alpha_M]$.
%%%%%%%%%%%%%%%%%%%%%%%%%%%%%%
%%%%%%%%%%%%%%%%%%%%%%%%%%%%%%
\subsection{Kontsevich's characteristic classes}\label{ss:kontsevich-class}
Let $\Gamma$ be a connected Jacobi diagram of degree $2n$ without a part like $\multimap$ and let $\omega(\Gamma)$ be the $3n(m-1)$-form on $C_{2n}(M)\gtimes\EDiff{M}$ defined by 
\[ \omega(\Gamma)\eqdef\bigwedge_{e:\,\mbox{\tiny edge of $\Gamma$}}\phi_e^*\alpha_{\Diff{M}} \]
where we fix a bijective correspondence between the set of vertices of $\Gamma$ and the set of $2n$ points in a configuration, and 
\[ \phi_e:C_{2n}(M)\gtimes\EDiff{M}\to C_2(M)\gtimes\EDiff{M}\]
is the projection in $C_{2n}(M)$-fibers corresponding to picking of the two endpoints of $e$. Note that the choice of the form $\alpha_{\Diff{M}}$ and therefore of $\omega(\Gamma)$ depends on the framing on $M$. Then the pushforward $(\pi_{C_{2n}(M)})_*\omega(\Gamma)$ along the fiber of $\pi_{C_{2n}(M)}$ yields an $n(m-3)$-form on $\BDiff{M}$. See Appendix~\ref{s:pushforward} for the definition of the pushforward.

According to \cite{Kon}, the form
\[ \zeta_{2n}\eqdef \sum_{\Gamma}\frac{(\pi_{C_{2n}(M)})_*\omega(\Gamma)[\Gamma]}{|\Aut{\Gamma}|}\in\Omega^{n(m-3)}(\BDiff{M};\calA_{2n}), \]
where the sum is over all connected trivalent graphs without $\multimap$ and $|\Aut\Gamma|$ be the order of the group of automorphisms of $\Gamma$, is closed and thus descends to an $\calA_{2n}$-valued universal characteristic class of framed smooth $M$-bundles. Further, $\R$-valued Kontsevich classes are defined by composing $\zeta_{2n}$ with any linear functional on $\calA_{2n}$. In the following, we will write $\zeta_{2n}(\sigma)=\int_{\sigma}\zeta_{2n}$ for a chain $\sigma$ on $\BDiff{M}$.

In the case $M$ is a 3-dimensional homology sphere, all the $\zeta_{2n}$ give rise to a universal $\R$-valued finite type invariants \cite{KT}.

%\clearpage

%%%%%%%%%%%%%%%%%%%%%%%%%%%%%%
%%%%%%%%%%%%%%%%%%%%%%%%%%%%%%
%%%%%%%%%%%%%%%%%%%%%%%%%%%%%%
\section{Bordism invariant of unframed $M$-bundles}\label{ss:bordism-inv}

In this section, we restrict our study mainly to smooth bundles with fiber a 5-dimensional homology sphere $M$. In this setting, we will show that the simplest Kontsevich class $\zeta_2$ after an addition of a certain multiple of the second signature defect invariant becomes a bordism invariant of unframed $M$-bundles. The strategy for the proof is mainly inspired by Lescop's nice explanation \cite{Les} of the Kuperberg-Thurston construction of unframed 3-manifold invariants and by Morita's construction of the secondary characteristic classes of surface bundles from the signature defects \cite{Mo}.

We restrict the holonomy group to the subgroup $\Diff' M\subset\mathrm{Diff}(M^\bullet\mbox{ rel }\partial)$ consisting of diffeomorphisms inducing homotopy trivial automorphisms on vertical framings. Namely, if $\varphi\in \mathrm{Diff}'M$, then $\varphi_*\tau_{M^\bullet}$ is homotopic to $\tau_{M^\bullet}$ for any vertical framing $\tau_{M^\bullet}$. This restriction does not lose the generality so much. Indeed, since $M^\bullet$ is a punctured homology sphere, the obstruction to homotopy two different framings on $M^\bullet$ lies in $H^5(M^\bullet,\partial M^\bullet;\pi_5(SO(5)))=H^5(M^\bullet,\partial M^\bullet;\Z_2)\cong\Z_2$. Hence $\varphi\circ\varphi$ for any $\varphi\in\Diff(M^\bullet\mbox{ rel }\partial)$ belongs to $\Diff'M$. We will call an $M$-bundle admitting a reduction of the holonomy to $\Diff'M$ a {\it $\Diff'M$-bundle}.

For a $\Diff'M$-bundle $\pi:E\to B$, we denote by $\pi^\bullet:E^\bullet\to B$ the bundle obtained from $E$ by restricting the fiber to $M^\bullet$. 
Let $\tau_{E^\bullet}$ denote a vertical framing of $E^\bullet$, if exists. Let $\pi_0:E_0\to B$ be the trivial $\Diff'M$-bundle $M\times B$ over $B$ vertically framed by the pullback of the framing on the fiber $E_{q_0}^\bullet= (\pi^{\bullet})^{-1}(q_0)$ by the projection $M\times B\to M\times q_0$ where $q_0\in B$ is the base point of $B$. Let $\overline{E}\eqdef E^\bullet\cup_{\partial=S^4\times B}(-E_0^\bullet)$ with vertical framing $\tau_{\overline{E}}=\tau_{E^\bullet}\cup_\partial(-\tau_{E_0^\bullet})$, which may be discontinuous at $\partial=S^4\times B$.\footnote{To avoid the discontinuity, one may instead take a vector bundle, which is not tangent to $\overline{E}$ near $\partial$. The resulting definition of the signature defect coincides with the previous one.} By Thom's cobordism theory, there exists a positive integer $N$ such that the disjoint union of $N$ copies of the 7-manifold $\overline{E}$ bounds a compact oriented $8$-manifold $W$, namely  $\partial W=\overline{E}\sqcup\cdots\sqcup \overline{E}$ ($N$ copies). 

Note that $TW|_{\overline{E}^{\sqcup N}}=(T\overline{E})^{\sqcup N}\oplus\varepsilon=(\pi^*TB\oplus\xi)^{\sqcup N}\oplus\varepsilon$ where $\xi$ is the vertical tangent bundle and $\varepsilon$ is the trivial 1-dimensional normal bundle over $\partial W=\overline{E}^{\sqcup N}$. Choose a connection on $TB$ and pull it back to $\pi^*TB$. Then together with the flat connection defined by $(\tau_{\overline{E}})^{\sqcup N}\oplus\tau_\varepsilon$, it defines a connection on $TW|_{\overline{E}^{\sqcup N}}$. Note that this flat connection is well-defined although $\tau_{\overline{E}}$ may be discontinuous since the vertical framings on both $E^\bullet$ and $E_0^\bullet$ are chosen so that they are standard near the boundaries and thus the induced flat connections near the boundaries are trivial. This connection can be extended to whole of $W$. The relative $L_2$-class is defined with this connection on $TW$ by Hirzebruch's $L$-polynomial given by
\[ 
	 L_2(TW;\tau_{\overline{E}}^{\sqcup N})=L_2(p_1,p_2)=\frac{1}{45}(7p_2-p_1^2)
 \]
where $p_j=p_j(TW;\tau_{\overline{E}}^{\sqcup N})$ is the $j$-th relative Pontrjagin class. It is known that the relative $p_2$-class can be interpreted as the obstruction class in $H^8(W,\partial W;\pi_7(SU(8)/SU(3)))$ to extend the partial 5-framing $\tau_{\overline{E}}^{\sqcup N}$ on $\partial W$ to the partial 5-framing over the complexified tangent bundle $TW\otimes\C$. In this setting, we assume that the sign convension for the relative $p_2$-class is determined by a fixed choice of the generator $[\mu]\in\pi_7(SU(8)/SU(3))\cong\Z$.
The second signature defect $\Delta_{2}(E;\tau_{E^\bullet})$ is defined by
\[ \Delta_{2}(E;\tau_{E^\bullet})\eqdef\frac{1}{N}\left[\int_W L_{2}(TW;\tau_{\overline{E}}^{\sqcup N})-\mathrm{sign}\,W\right].\]
The $k$-th signature defect $\Delta_k$ for $k>2$ is also defined by using $L_k$.
%\clearpage

\begin{Prop}\label{prop:L2defined}
$\Delta_{2}(E;\tau_{E^\bullet})$ is well-defined. That is $\Delta_{2}(E;\tau_{E^\bullet})$ is independent of the choices of the connection, the bounding manifold $W$ and the number $N$ of copies.
\end{Prop}
The proof of Proposition~\ref{prop:L2defined} is the same as \cite[Proposition~7.3]{Mo}.

\begin{Thm}\label{thm:A}Any $\Diff'M$-bundle over a closed connected oriented 2-manifold with fiber a 5-dimensional homology sphere, can be vertically framed. Moreover in this case, the number
\[ \hat{\zeta}_2(E)\eqdef \zeta_2(E;\tau_{E^\bullet})-\frac{15}{112}\,\Delta_2(E;\tau_{E^{\bullet}})\,[\Theta]\in\calA_2\]
does not depend on the choice of a vertical framing $\tau_{E^\bullet}$ and is a bordism invariant $\Omega_2(B\mathrm{Diff}'{M})\to\calA_2$ of smooth unframed $\Diff'M$-bundles.
\end{Thm}

Theorem~\ref{thm:A} immediately implies
\begin{Cor}
For $\dim{M}=5$, $\zeta_2$ is non-trivial.
\end{Cor}

We do not know whether $B\Diff'M$ has the homotopy type of a finite CW-complex (while for 3-dimensional manifolds, it is known to be true, which was conjectured by Kontsevich and proved by Hatcher and McCullough \cite{HM}). So we do not know whether $\hat{\zeta}_2$ descends to a cohomology class. 

By a similar argument as in \cite{KT, Les}, we have
\begin{equation}\label{eq:dzeta}
 d\zeta_{2n}=\sum_{\Gamma}\frac{[\Gamma]}{|\Aut\Gamma|}\int_{S_{2n}(TM)_b}\omega(\Gamma)\quad (b\in\BDiff{M}), 
\end{equation}
which vanishes on $\BDiff{M}$. Here $S_{2n}(TM)_b\to M_b$ denotes the bundle associated to $TM_b$ whose fiber is the space of configurations of $2n$ points in a $5$-dimensional plane modulo overall translations and dilations. Indeed, $d\zeta_{2n}$ evaluated on any $(2n+1)$-chain $\sigma$ can be expressed as an integral of a pulled back form from the fiber of a point of $\sigma$. Then the integral vanishes by a dimensional reason.  

If one wants to make $\zeta_{2n}$ framing independent, it suffices to add some correction term to cancel the RHS of (\ref{eq:dzeta}). Theorem~\ref{thm:A} says that $\frac{15}{112}\,\Delta_2(E;\tau_{E^\bullet})$ is a suitable correction.

\begin{Rem}\label{rem:triviality}
We do not know whether $\hat{\zeta}_2$ of Theorem~\ref{thm:A} is non-trivial or not when $M$ is 5-dimensional. If it is trivial, then one obtains the relation 
\[ \zeta_2(E;\tau_{E^\bullet})=\frac{15}{112}\Delta_2(E;\tau_{E^\bullet})[\Theta]. \]
If it is non-trivial, then it is expected that it measures independent structures of bundles that does not determined by the signature defect. 

When $M$ is $(2k-1)$-dimensional at least 7, unframed invariant $\hat{\zeta}_2$ is also obtained in Theorem~\ref{thm:AA} for bundles over $S^{2(k-2)}$. In particular, when $2k-1=7$, we obtain the result about the non-triviality of $\hat{\zeta}_2$ (Corollary~\ref{cor:z2-nontriv}).

The formula for $\hat{\zeta}_2$ is similar to Morita's splitting formula for the Casson invariant \cite{Mo}:
\[ \lambda(M)=\frac{1}{6}\int_{C_2(M)}\alpha_M^3-\frac{1}{24}\cdot 3\Delta_1(M;\tau_M). \]
(This is the version described in \cite{KT,Les2}). So the existence of analogous properties and constructions for the Casson invariant may be expected for $\hat{\zeta}_2$.\qed
\end{Rem}

%\clearpage

\subsection{Vertical framings for $\Diff'M$-bundles}
The following proposition proves the first part of Theorem~\ref{thm:A} and allows us to define the Kontsevich classes for any $\Diff'M$-bundle over a closed connected oriented 2-manifold.
\begin{Prop}\label{prop:E-framed}
Let $M$ be a $5$-dimensional homology sphere. Any $\Diff'M$-bundle over a closed connected oriented $2$-manifold can be vertically framed.
\end{Prop}
\begin{proof}
Let $\pi:E\to B$ be a $\Diff'M$-bundle over a closed connected oriented 2-manifold $B$ and choose a cell decomposition of $B$ with one 0-cell. Since the holonomy is contained in $\Diff'M$, the vertical framing extends over the 1-skeleton.

To see that the vertical framing extends over a 2-skeleton of $B$, we consider a trivial $\Diff'M$-bundle $M\times e^2 \to e^2$ over the 2-cell $e^2$ and consider the obstruction for homotoping the standard vertical framing over $C=\partial e^2\cong S^1$ into the vertical framing over the image of $C$ under the attaching map, which is determined by the above extension over the 1-skeleton. We can choose a vertical framing of the trivial bundle over $e^2$ so that the two vertical framings coincide at the fiber over the base point $q_0$ of $C$. The difference of the two vertical framings $\tau_1, \tau_2$ can be considered as a map
\[ g=\tau_2\circ\tau_1^{-1}:M^\bullet\times C\to GL_+(\R^5) \]
which is trivial on $(M^\bullet\times \{q_0\})\cup(\partial M^\bullet\times C)$. Moreover, this map is reduced modulo homotopy to a map into $SO(5)\subset GL_+(\R^5)$ by the deformation retraction given by the Gram-Schmidt orthonormalization. It suffices to prove the vanishing of the obstruction for homotoping $g$ into the constant map.

Choose a cell decomposition of $M^\bullet\times C$ with respect to its boundary, induced by a cell decomposition of $M^\bullet$ with respect to the boundary. By Lemma~\ref{lem:H(E,dE)} below, we have
\[ H^j(M^\bullet\times C, (M^\bullet\times\{q_0\})\cup(\partial M^\bullet\times C);\pi_j(SO(5)))=0\]
for $2\leq j\leq 6$, which implies that the homotopy extends over the whole of $M^\bullet\times C$. Namely, the vertical framing extends over whole of $B$.
\end{proof}

\begin{Lem}\label{lem:H(E,dE)}Let $\pi:E\to B$ be a $\Diff'M$-bundle over a closed connected oriented manifold $B$ of dimension $\leq 2$. Then
\[ H^i(E^\bullet,\partial{E}^\bullet\cup E_{q_0}^\bullet;\pi_i(SO(5)))=0 \]
for $0\leq i\leq 6$.
\end{Lem}
\begin{proof}
First we compute the homology group $H_i(E^\bullet,\partial{E}^\bullet\cup E_{q_0}^\bullet;\Z)$ via the homology exact sequence
\begin{equation}\label{eq:homology-seq}
 \mapsright{}
	H_i(\partial{E}^\bullet\cup E_{q_0}^\bullet;\Z)\mapsright{}
	H_i(E^\bullet;\Z)\mapsright{}
	H_i(E^\bullet,\partial{E}^\bullet\cup E_{q_0}^\bullet;\Z)\mapsright{}. 
\end{equation}
Since the $\Diff'M$-bundle $\pi^\bullet:E^\bullet\to B$ is homologically a disk bundle, its homology is isomorphic to that of $B$:
\begin{equation}\label{eq:H(E)}
 H_i(E^\bullet;\Z)\cong H_i(B;\Z). 
\end{equation}
The homology of $\partial{E}^\bullet\cup E_{q_0}^\bullet$ is determined via the Mayer-Vietoris sequence as
\begin{equation}\label{eq:H(dE)}
 H_i(\partial{E}^\bullet\cup E_{q_0}^\bullet;\Z)\cong\left\{
	\begin{array}{ll}
		H_i(B;\Z) & \mbox{if }0\leq i\leq 2\\
		H_{i-4}(B;\Z) & \mbox{if }5\leq i\leq 6\\
		0 & \mbox{otherwise}
	\end{array}
	\right. 
\end{equation}
Substituting (\ref{eq:H(E)}) and (\ref{eq:H(dE)}) into (\ref{eq:homology-seq}), we have $H_i(E^\bullet,\partial{E}^\bullet\cup E_{q_0}^\bullet;\Z)=0$ for $0\leq i\leq 5$. Hence by the universal coefficient theorem, we have
\[ H^i(E^\bullet,\partial{E}^\bullet\cup E_{q_0}^\bullet;\pi_i(SO(5)))=0 \]
for $0\leq i\leq 5$. Furthermore, by $\pi_6(SO(5))=0$, we have
\[ H^6(E^\bullet,\partial{E}^\bullet\cup E_{q_0}^\bullet;\pi_6(SO(5)))=0. \]
\end{proof}

%%%%%%%%%%%%%%%%%%%%%%%%%%%%%%
%%%%%%%%%%%%%%%%%%%%%%%%%%%%%%
\subsection{Framing dependence of $\zeta_2$}

The proof of Theorem~\ref{thm:A} is carried out by showing that the framing dependences of $\zeta_2$ and of $\Delta_2$ differ by a constant multiple. In this subsection we compute the difference of $\zeta_2$ for two different vertical framings. 

\begin{Lem}\label{lem:z-z=dz}Let $(\pi:E\to B,\tau_{E^\bullet})$ a vertically framed null bordant $\Diff'M$-bundle over a not necessarily connected closed 2-manifold $B$. Then $\zeta_2(E;\tau_{E^\bullet})=0$.
\end{Lem}
\begin{proof}
Since $\zeta_2$ is a cocycle on $\BDiff{M}$, it is a framed bordism invariant. Thus the result follows.
\end{proof}

\begin{Lem}\label{lem:dependsonhtpy}
$\zeta_2(E;\tau_{E^\bullet})$ depends only on the homotopy class of $\tau_{E^\bullet}$.
\end{Lem}
\begin{proof}
Let $\tau_{E^\bullet}$ and $\tau_{E^\bullet}'$ be two mutually homotopic vertical framings. We prove $\zeta_2(E;\tau_{E^\bullet})=\zeta_2(E;\tau_{E^\bullet}')$.

The homotopy gives rise to a cylinder $E\times I$ with a vertical framing $\tilde{\tau}_{E^\bullet}(t)\ (t\in I)$ such that $\tilde{\tau}_{E^\bullet}(0)=\tau_{E^\bullet}$ and $\tilde{\tau}_{E^\bullet}(1)=\tau_{E^\bullet}'$. Existence of the framed cylinder $E\times I$ implies that the two vertically framed bundles $(E;\tau_{E^\bullet})$ and $(E;\tau'_{E^\bullet})$ are vertically framed bordant. Hence Lemma~\ref{lem:z-z=dz} concludes the proof.
\end{proof}

%\clearpage

\begin{Lem}\label{cor:homotopy-obstruction}
Let $\pi:E\to B$ denote a $\Diff'M$-bundle over a closed connected oriented $2$-manifold $B$. Then there is a homotopy deforming any continuous map
\[ g:E^\bullet\to SO(5) \]
which is trivial on $\partial{E}^\bullet\cup E_{q_0}^\bullet$ into a map that is trivial outside a $7$-ball embedded into $E^\bullet$.
\end{Lem}
\begin{proof}
Lemma~\ref{lem:H(E,dE)} implies that the homotopy extends from $\partial{E}^\bullet\cup E_{q_0}^\bullet$ over the $6$-skeleton of $E^\bullet$.
\end{proof}

For any map $G: (E^\bullet, \partial E^\bullet\cup E_{q_0}^\bullet)\to (SO(5),1)$, let $\psi(G): \R^5\times E^\bullet\to \R^5\times E^\bullet$ be defined by $\psi(G)(v,x)\eqdef (G(x)v,x)$. 

\begin{Lem}\label{lem:zeta-zeta}
Let $\pi:E\to B$ be a $\Diff'M$-bundle over a closed connected oriented 2-manifold $B$ and let $\tau_{E^\bullet}$ be a vertical framing on it. Then $\zeta_2(E;\psi(G)\circ \tau_{E^\bullet})-\zeta_2(E;\tau_{E^\bullet})$ does not depend on $\tau_{E^\bullet}$. It depends only on the homotopy class of $\psi(G)$.
\end{Lem}
\begin{proof}
Let $\tilde{\pi}:\widetilde{E}\to\widetilde{B}$ be the trivial $\Diff'M$-bundle over the cylinder $\widetilde{B}=B\times I$ so that $\widetilde{E}=E\times I$. Suppose that, $E\times \{1\}$ and $E\times \{0\}$ in $\widetilde{E}$ are vertically framed by $\psi(G)\circ\tau_{E^\bullet}$ and $\tau_{E^\bullet}$ respectively. 

By Lemma~\ref{cor:homotopy-obstruction}, we may assume after a homotopy that $\psi(G)\circ\tau_{E^\bullet}$ and $\tau_{E^\bullet}$ coincide outside $\pi^{-1}B^{2}$ where $B^{2}\subset B$ is an embedded $2$-disk. In other words, the vertical framing over $\partial \widetilde{B}=B\times\{1\}\sqcup-B\times\{0\}$ extends to $\widetilde{B}$ outside an embedded $3$-ball $B^3\subset \widetilde{B}$. Furthermore we may consider $\widetilde{E}^\circ\eqdef\tilde{\pi}^{-1}{\widetilde{B}\setminus \mathrm{int}(B^{3})}$ as a cobordism between $E\sqcup (M\times S^2)$ and $-E$ vertically framed by $\tau_{\widetilde{E}^{\circ\bullet}}$. We denote by $\tau_G$ the induced vertical framing on $M\times S^{2}$.
%\begin{figure}
%\psfrag{B2}[cc][cc]{$B^2$}
%\psfrag{B3}[cc][cc]{$B^3$}
%\psfrag{B}[cc][cc]{$B$}
%\psfrag{-B}[cc][cc]{$-B$}
%\psfrag{S2}[cc][cc]{$S^2$}
%\fig{eus-e.eps}
%\caption{}\label{fig:eus-e}
%\end{figure}

By Lemma~\ref{lem:z-z=dz}, we have
$	 \zeta_2(E;\psi(G)\circ\tau_{E^\bullet})+\zeta_2(M\times S^{2};\tau_G)-\zeta_2(E;\tau_{E^\bullet})
	=0$. Namely, $\zeta_2(E;\psi(G)\circ\tau_{E^\bullet})-\zeta_2(E;\tau_{E^\bullet})=-\zeta_2(M\times S^{2};\tau_G)$ does not depend on $\tau_{E^\bullet}$.
\end{proof}

%\clearpage

The last proposition allows us to define
\[ \zeta'_2(E;G)\eqdef \zeta_2(E;\psi(G)\circ \tau_{E^\bullet})-\zeta_2(E;\tau_{E^\bullet}). \]

Let $p:E_{\rho}\to S^{8}$ be the real $5$-dimensional vector bundle over $S^{8}=B^{8}\cup_{\partial=S^{7}}(-B^{8})$ defined by
\[ E_{\rho}\eqdef ( \R^{5}\times B^{8})\cup_h(\R^{5}\times -B^{8}) \]
where the gluing map $h:\R^{5}\times \partial B^{8}=\R^5\times S^{7}\to \R^5\times S^{7}$ is the twist defined by $(v,x)\mapsto (\rho(x)v,x)$ with a smooth map $\rho:S^{7}\to SO(5)\subset GL_+(\R^{5})$ representing the generator of $\pi_{7}(SO(5))\cong\Z$. Here we choose $\rho$ so that $[\rho]\in\pi_7(SO(5))$ is mapped by the inclusion $i_*:\pi_7(SO(5))\to\pi_7(SU(8)/SU(3))=\langle[\mu]\rangle$ to a positive multiple of $[\mu]$. Non-triviality of $i_*[\rho]$ will be shown later.

For an $\R^5$ vector bundle $E$, we denote by $S_2(E)$ the $S^4$-bundle associated to $E$. Let $\omega_T$ be a closed $4$-form on the $S^{4}$-bundle $S_2(E_{\rho})$ representing the Thom class such that $\iota^*\omega_T=-\omega_T$ under the involution $\iota:E_{\rho}\to E_{\rho}$ defined by $\iota(x,v)=(x,-v)$. Let
\[ 	\delta_{2}(E_{\rho})\eqdef 
		\frac{[\Theta]}{|\Aut\Theta|}\int_{S_{2}(E_{\rho})}\omega_T^3
		=\frac{[\Theta]}{12}\int_{S_{2}(E_{\rho})}\omega_T^3\in\calA_{2}.\]
One can prove that $\delta_2(E_{\rho})$ does not depend on the choice of $\omega_T$ satisfying $\iota^*\omega_T=-\omega_T$ within the cohomology class because $S_2(E_\rho)$ is a closed manifold.

\begin{Lem}\label{lem:zeta(rho)}
$\zeta'_2(E;G)=\delta_2(E_{\rho})$ if $G$ is homotopic to a map $G_E(\rho)$ that coincides with $\mathrm{id}$ outside some embedded $7$-ball $B^{7}$ in $E$ and the image of $B^{7}$ under $G_E(\rho)$ is homotopic to $\rho$.
\end{Lem}
\begin{proof}
By Stokes' theorem and by (\ref{eq:dzeta}), we have
\[ \begin{split}
	\zeta_2'(E;G_E(\rho))&=\zeta_2(E;\psi(G_E(\rho))\circ \tau_{E^\bullet})-\zeta_2(E;\tau_{E^\bullet})\\
	&=\frac{[\Theta]}{12}\int_{S_2(\R^5\times I\times B^7)}\omega_{T}(\R^5\times I\times B^7)^3
	\end{split} 
\]
where $\omega_T(\R^5\times I\times B^7)$ denotes the 4-form representing the Thom class of the associated $S^4$-bundle $S^4\times I\times B^7$ extending $\omega(\tau_{E^\bullet})$ and $\omega(\psi(G_E(\rho))\circ\tau_{E^\bullet})$ given on $S^4\times B^7\times \partial I$ such that the involution $\iota^*$ acts as $-1$. Existence of such a 4-form is because the restriction induces an isomorphism from $H^4(S^4\times I\times B^7;\R)$ to $H^4(\partial(S^4\times I\times B^7);\R)$. 

On the other hand, $\delta_2(E_{\rho})$ can be computed from the definition as follows. Consider the base $S^8$ of $S_2(E_\rho)$ as the union $B^8\cup_{\partial=S^7=B^7\cup_\partial(-B^7)}(S^7\times I)\cup_\partial (-B^8)$. We can give a trivialization $\tau_{B^8}$ of $S_2(E_\rho)$ over $S^8\setminus \mathrm{int}(I\times B^7)$, and we can choose $\omega_T$ so that it is an extension of the spherical form determined by $\tau_{B^8}$ on $S^8\setminus \mathrm{int}(I\times B^7)$. Then we have
\[ \begin{split}
	\delta_2(E_{\rho})&=\frac{[\Theta]}{12}\int_{S_2(E_{\rho})}\omega_T^3
	=\frac{[\Theta]}{12}\Bigl(\int_{S_2(\R^5\times I\times B^7)}\omega_{T}^3
	+\int_{S_2( \R^5\times (S^8\setminus\mathrm{int}(I\times B^7)))}\omega_{T}^3\Bigr)\\
	&=\frac{[\Theta]}{12}\int_{S_2(\R^5\times I\times B^7)}\omega_{T}^3
		=\frac{[\Theta]}{12}\int_{S_2(\R^5\times I\times B^7)}\omega_{T}(\R^5\times I\times B^7)^3\\
	&	=\zeta_2'(E;G_E(\rho))
	\end{split} \]
where the third equality follows from a dimensional reason and the fourth equality follows from the fact that the $\omega_T$ on $S_2(E_{\rho})$ can be chosen as an extension of $\omega_T(\R^5\times I\times B^7)$.
\end{proof}

%\clearpage

For a $\Diff'M$-bundle $\pi:E\to B$, we denote by $[E,SO(5)]^\bullet$ the set of homotopy classes of continuous maps 
\[G:(E^\bullet,\partial{E}^\bullet\cup E_{q_0}^\bullet)\to (SO(5),1).\]
The following proposition is a key to prove Theorem~\ref{thm:A}, describing the structure of the set of homotopy classes of vertical framings.

\begin{Prop}\label{prop:[E,SO(5)]}Let $\pi:E\to B$ be a vertically framed $\Diff'M$-bundle over a closed connected oriented $2$-manifold. Then $[E,SO(5)]^\bullet=\langle [G_E(\rho)]\rangle$ ($G_E(\rho)$ is defined in Lemma~\ref{lem:zeta(rho)}), the free abelian group generated by $G_E(\rho)$. Thus the degree in $[E,SO(5)]^\bullet$ is defined by $p[G_E(\rho)]\mapsto p$.
\end{Prop}
\begin{proof}
By Lemma~\ref{cor:homotopy-obstruction}, the obstruction to homotoping $G$ into the constant map over whole of $E$ is described by a homotopy class of a map $\partial (B^{7}\times I)\cong S^{7}\to SO(5)$, which is an element of $\pi_{7}(SO(5))\cong\Z$. 
\end{proof}

%\clearpage

\begin{Lem}\label{lem:dzeta=d2}
Let $G\in [E,SO(5)]^\bullet$. Then we have
\[ \zeta_2'(E;G)=\delta_2(E_{\rho})\,\deg G. \]
\end{Lem}
\begin{proof}
By Lemma~\ref{lem:zeta-zeta}, we have
\[ \begin{split}
	\zeta_2'(g)+\zeta_2'(h)
	&=(\zeta_2(E;\psi(g)\circ\psi(h)\circ\tau_{E^\bullet})
		-\zeta_2(E;\psi(h)\circ\tau_{E^\bullet}))\\
	&\quad +(\zeta_2(E;\psi(h)\circ\tau_{E^\bullet})
		-\zeta_2(E;\tau_{E^\bullet}))=\zeta_2'(E;gh).
	\end{split}\]
Therefore $\zeta_2':[E,SO(5)]^\bullet\to\calA_2$ is a group homomorphism. Then by Proposition~\ref{prop:[E,SO(5)]}, $\zeta_2'$ is a multiple of $\deg$ with some constant in $\calA_2$. Lemma~\ref{lem:zeta(rho)} implies that the constant is exactly equal to $\delta_2(E_{\rho})$.
\end{proof}

%\clearpage

%%%%%%%%%%%%%%%%%%%%%%%%%%%%%%
%%%%%%%%%%%%%%%%%%%%%%%%%%%%%%
\subsection{Framing dependence of Pontrjagin numbers}

As for $\zeta_2$, we compute the difference between the relative Pontrjagin numbers for two different vertical framings. We only need to see the framing dependence of the second Pontrjagin number because the relative cohomology $H^4(E^\bullet \times I,\partial(E^\bullet \times I);\Z)$ vanishes and the difference of square of the first relative Pontrjagin number vanishes. 

\begin{Lem}\label{lem:p-p}Let $\pi:E\to B$ is a $\Diff'M$-bundle over a closed connected oriented $2$-manifold. Then $p_2(E;\psi(G)\circ \tau_{E^\bullet})-p_2(E;\tau_{E^\bullet})$ does not depend on $\tau_{E^\bullet}$. It depends only on the homotopy class of $\psi(G)$.
\end{Lem}
\begin{proof}
The difference computes the second relative Pontrjagin number of $E^\bullet\times I$ with respect to the vertical framings $\psi(G)\circ\tau_{E^\bullet}$ and $\tau_{E^\bullet}$ given on $E^\bullet\times\{0,1\}$, and with respect to the standard vertical framing on $\partial E^\bullet\times I$. Then the proof may be carried out by a similar argument as in Lemma~\ref{lem:zeta-zeta} with the fact that the second relative Pontrjagin number vanishes on vertically framed cobordisms.
\end{proof}
%\clearpage

Lemma~\ref{lem:p-p} allows us to define
\[ p'_2(E;G)\eqdef p_2(E;\psi(G)\circ \tau_{E^\bullet})-p_2(E;\tau_{E^\bullet}). \]
\begin{Lem}\label{lem:p'2}Let $\pi:E\to B$ be $\Diff'M$-bundle over a closed connected oriented 2-manifold. Then
\begin{equation}\label{eq:p'=48d(G)} p'_2(E;G)=48\,\deg G. 
\end{equation}
\end{Lem}
\begin{proof}
Since $p_2'(E;G):[E,SO(5)]^\bullet\to\Q$ is a group homomorphism, it follows from Proposition~\ref{prop:[E,SO(5)]} that 
\[ p_2'(E;G)=p_2'(E;G_E(\rho))\deg{G}. \]
So it suffices to prove that $p_2'(E;G_E(\rho))=48$.

The second relative Pontrjagin class $p_2'$ is considered as the obstruction to extend the vertical framing on $\partial(E^\bullet\times I)$ to the complexified vertical tangent bundle of $E^\bullet\times I$. This obstruction lies in $H^7(E^\bullet,\partial E^\bullet;\pi_7(SU(8)/SU(3)))=H^7(E^\bullet,\partial E^\bullet; \Z)$. In the case of $G_E(\rho)$, the obstruction is the image of $[\rho]\in\pi_7(SO(5))$ under the inclusion $\pi_7(SO(5))\to\pi_7(SU(8)/SU(3))$. This inclusion factors through $\pi_7(SU(5))\cong\Z$ and the following two lemmas conclude the proof.
\end{proof}

\begin{Lem}\label{lem:SU->V}
The natural inclusion $i:SU(5)\to SU(8)/SU(3)$ sends the generator of \\
$\pi_7(SU(5))\cong\Z$ to $\pm 6$ times the generator of $\pi_7(SU(8)/SU(3))\cong\Z$.
\end{Lem}
\begin{proof}
This is a direct consequence of the following homotopy sequence of the bundle:
\[ \begin{array}{cccccccc}
 \mapsright{} & \pi_7(SU(5)) & \mapsright{i_*} & \pi_7(SU(8)/SU(3)) & \mapsright{} & \pi_7(SU(8)/(SU(5)\times SU(3))) & \mapsright{} & 0 \\
	& \rotatebox[origin=c]{90}{$\cong$} & & \rotatebox[origin=c]{90}{$\cong$} & & \rotatebox[origin=c]{90}{$\cong$} && \\
	& \Z & & \Z & & \Z_6 & &
\end{array}\]
\end{proof}

The following lemma follows from a result in \cite{Lun}.
\begin{Lem}
The natural inclusion $c:SO(5)\to SU(5)$ sends the generator of $\pi_7(SO(5))\cong \Z$ to $\pm 8$ times the generator of $\pi_7(SU(5))\cong\Z$.
\end{Lem}

%\clearpage

%%%%%%%%%%%%%%%%%%%%%%%%%%%%%%
%%%%%%%%%%%%%%%%%%%%%%%%%%%%%%
\subsection{Computation of $\delta_2(E_{\rho})$ and framing correction}

The following lemma is proved in \cite{BC}.
\begin{Lem}[Bott-Cattaneo]\label{lem:bott-cattaneo}
Let $\pi:E\to B$ be an $\R^{2k-1}$-vector bundle and let $S(E)$ be its associated sphere bundle with $e\in H^{2k-2}(S(E);\R)$ be the canonical Euler class. Then
\[ \pi_* e^3=2p_{k-1}(E), \]
twice of the $(k-1)$-st Pontrjagin class.
\end{Lem}

Since the Euler number of $S^{2k}$ is 2, $e$ restricts to twice the generator of $H^{2k}(S^{2k};\Z)$.

\begin{Lem}\label{lem:d=theta}$\delta_2(E_{\rho})=[\Theta]. $
\end{Lem}
\begin{proof}
The integral part of $\delta_2(E_{\rho})$: $\int_{S_2(E_{\rho})}\omega^3$ is equal to $\frac{2}{2^3}\langle p_2(S_2(E_\rho)),[S^8]\rangle=\frac{1}{4}\langle p_2(S_2(E_\rho)),[S^8]\rangle$ by Lemma~\ref{lem:bott-cattaneo}.

We can choose a partial 5-frame outside an 8-ball $B^8$ embedded into the base $S^8$ of $E_{\rho}$. By the construction of $E_{\rho}$, the obstruction class to extend the partial 5-frame over whole of $S^8$, which lies in $H^8(S^8,B^8;\pi_7(SO(5)))$, is defined as the map sending the boundary $\partial B^8$ of the only 8-cell into the image of $\rho:S^7\to SO(5)$. Recall that the second Pontrjagin class can be considered as this obstruction class, which can be considered lying in $H^8(S^8,B^8;\pi_7(SU(8)/SU(3))=\Z)$.

According to Lemma~\ref{lem:p'2}, $[\rho]\in\pi_7(SO(5))$ is mapped under the inclusion $\pi_7(SO(5))\to\pi_7(SU(8)/SU(3))$ into $48$ times the generator of $\pi_7(SU(8)/SU(3))$. Hence $\langle p_2(S_2(E_\rho)),[S^8]\rangle=48$ and
\[ \delta_2(E_{\rho})=\frac{[\Theta]}{12}\int_{S_2(E_{\rho})}\omega^3=\frac{[\Theta]}{12}\cdot \frac{48}{4} =[\Theta]. \]
\end{proof}
%\clearpage
%%%%%%%%%%%%%%%%%%%%%%%%%%%%%%
%%%%%%%%%%%%%%%%%%%%%%%%%%%%%%

Now we shall see in the case of bundles with fiber a 5-dimensional homology sphere over a base closed 2-dimensional manifold, that unframed bordant implies vertically framed bordant.

\begin{Lem}\label{lem:bordism-framed}
If two vertically framed $\Diff'M$-bundles $\pi_j:E_j\to B_j\ (j=0,1)$ over closed oriented $2$-dimensional manifolds $B_j$ are unframed bordant, i.e., they define the same element of $\Omega_2(B\Diff'{M})$, then there exists a $\Diff'M$-bundle $\tilde{\pi}: \widetilde{E}\to\widetilde{B}$ such that 
\begin{enumerate}
\item $\partial \widetilde{B}=B_1\sqcup (-B_0)$,
\item $\tilde{\pi}|_{\partial \widetilde{E}}=\pi_1\sqcup (-\pi_0)$,
\end{enumerate}
and any given vertical framing over $\partial \widetilde{B}$ that is standard at the base point extends to a vertical framing over $\widetilde{B}\setminus B^{3}$ for some embedded $3$-disk $B^{3}\subset\widetilde{B}$.
\end{Lem}
\begin{proof}
Existence of $\tilde{\pi}:\widetilde{E}\to\widetilde{B}$ satisfying (1) and (2) is clear. For the last assertion, choose a cell decomposition of $\widetilde{B}$ with respect to $\partial\widetilde{B}$. The same argument as in Proposition~\ref{prop:E-framed} shows that the vertical framing also extends over the $2$-skeleton of $\widetilde{B}$.
\end{proof}

\begin{Cor}\label{cor:bordism-framed}
If two $\Diff'M$-bundles $\pi_j:E_j\to B_j\ (j=0,1)$ over closed oriented $2$-dimensional manifolds $B_j$ are unframed bordant, i.e., they define the same element of $\Omega_2(B\Diff'{M})$, then there exists a $\Diff'M$-bundle $\tilde{\pi}: \widetilde{E}\to\widetilde{B}$ and a vertical framing $\tau_{\tilde{\pi}^{-1}\partial \widetilde{B}}$ over $\partial\widetilde{B}$ such that
\begin{enumerate}
\item $\partial \widetilde{B}=B_1\sqcup (-B_0)$,
\item $\tilde{\pi}|_{\partial \widetilde{E}}=\pi_1\sqcup (-\pi_0)$,
\item $\tau_{\tilde{\pi}^{-1}\partial \widetilde{B}}$ extends over $\widetilde{B}$.
\end{enumerate}
\end{Cor}
\begin{proof}
Choose any vertical framing over $\partial\widetilde{B}$. Then the vertical framing extends to $\widetilde{B}\setminus B^{3}$ by Lemma~\ref{lem:bordism-framed}. After a homotopy, we may assume that the trivial bundle $\tilde{\pi}^{-1}B^{3}$, where the obstruction may be included, lies in a thin cylinder $\tilde{\pi}^{-1}(B_1\times [1-\varepsilon,1])$ near $B_1$. Then cut off the cylinder $\tilde{\pi}^{-1}(B_1\times (1-\varepsilon,1])$ from $\tilde{\pi}$. The resulting bundle is the desired vertically framed bordism.
\end{proof}

%\clearpage

Since $\zeta_2$ is a cocycle on $\BDiff M$, it is a vertically framed bordism invariant. Further we can also prove the following
\begin{Prop}
$\Delta_2$ is a vertically framed bordism invariant of $\Diff'M$-bundle over a closed oriented 2-manifold.
\end{Prop}
\begin{proof}
Let $\pi_j:E_j\to B_j\ (j=0,1)$, $\tilde{\pi}:\widetilde{E}\to\widetilde{B}$ and the vertical framing $\tau_{\widetilde{E}^\bullet}$ on $\widetilde{E}^\bullet$ be as in Corollary~\ref{cor:bordism-framed}. Note that any connections on $TB_j$ can be extended over $T\widetilde{B}$. We show that the signature defect $\Delta_{2}$ for $W=\widetilde{E}^\bullet$ vanishes.

We have $\mathrm{sign}\,\widetilde{E}^\bullet=0$ because $H^{4}(\widetilde{E}^\bullet;\Q)=0$. The first relative Pontrjagin class $p_1(T'\widetilde{E}^\bullet;\tau_{\widetilde{E}^\bullet})$ also vanishes because $H^4(\widetilde{E}^\bullet,\partial\widetilde{E}^\bullet;\Z)=0$. Further, $p_2(T'\widetilde{E}^\bullet;\tau_{\widetilde{E}^\bullet})=0$ because $\widetilde{E}^\bullet$ is vertically framed.
\end{proof}
%\clearpage

\begin{proof}[Proof of Theorem~\ref{thm:A}]
By Lemma~\ref{lem:dzeta=d2}, \ref{lem:p'2}, the vertically framed bordism invariant $\hat{\zeta}_2:\Omega_2(\BDiff{M})\to\calA_2$ defined by
\[ \begin{split}
	\hat{\zeta}_2(E)&=\zeta_2(E;\tau_{E^\bullet})
	-\frac{45}{7}\cdot\frac{1}{48}\Delta_2(E;\tau_{E^\bullet})\delta_2(E_{\rho})\\
	&=\zeta_2(E;\tau_{E^\bullet})-\frac{15}{112}\Delta_2(E;\tau_{E^\bullet})[\Theta]
	\end{split} \]
does not depend on the vertical framing $\tau_{E^\bullet}$. Corollary~\ref{cor:bordism-framed} says that unframed bordant implies vertically framed bordant. Thus $\hat{\zeta}_2$ can be considered as an unframed bordism invariant $\Omega_2(B\Diff'M)\to\calA_2$. Note that in the case $N>1$ in the definition of the signature defect, the $\deg{G}$ in (\ref{eq:p'=48d(G)}) may become $N$ times as much as the connected case. Then this cancels with the $\frac{1}{N}$ factor in the definition of the signature defect.
\end{proof}

%%%%%%%%%%%%%%%%%%%%%%%%%%%%%%
\subsection{Unframed homotopy invariant for $M$-bundle with higher $\dim{M}$}

We consider any $M$-bundle over $S^i$ coming from an element of $\pi_i B\Diff(M^\bullet\mbox{ rel }\partial)$. Thus sums or integer multiples of bundles are defined.

\begin{Thm}\label{thm:AA}
Let $M$ be a $(2k-1)$-dimensional homology sphere with $k\geq 4$ and let $\pi:E\to S^{2(k-2)}$ be an $M$-bundle over $S^{2(k-2)}$. Then there exists unique positive integer $p_k$ such that $p_k\pi:p_kE\to S^{2(k-2)}$ can be vertically framed for all $\pi$. Further, if $\tau_{E^\bullet}$ is a vertical framing on $p_k\pi$, then the number
\[ \hat{\zeta_2}(E)\eqdef
	\zeta_2(p_kE; \tau_{E^\bullet})
	-\frac{(2k-2)!}{3\cdot 2^{2k+2}(2^{2k-3}-1)B_{k-1}}\Delta_{k-1}(p_kE;\tau_{E^\bullet})[\Theta]
	\in\calA_2, \]
where $B_{k-1}$ is the $(k-1)$-th Bernoulli number, does not depend on the choices of $\tau_{E^\bullet}$, and is a homotopy invariant $\pi_{2(k-2)}B\Diff(M^\bullet\mbox{ rel }\partial)\otimes \R\to\calA_2$ of unframed $M$-bundles.
\end{Thm}

\begin{proof}[Outline of Proof]
By a similar argument as in Lemma~\ref{lem:H(E,dE)}, we can show that 
\[ H_i(E^\bullet,\partial E^\bullet\cup E^\bullet_{q_0};\Z)\cong
	\left\{\begin{array}{ll}
		\Z & \mbox{if $i=4k-5(=2(k-2)+(2k-1))$}\\
		0 & \mbox{otherwise}
	\end{array}\right. \]
(The assumption that the base is $S^{2(k-2)}$ is used here to simplify the homology.) This together with $\pi_{4k-6}(SO(2k-1))\otimes\Q=0$ implies that there exists unique $p_k$ such that $p_k\pi$ can be vertically framed for all $\pi$. 

That $\hat{\zeta_2}(E)$ does not depend on $\tau_{E^\bullet}$ follows by the same argument as in the previous subsection with the following facts:
\begin{itemize}
\item $[p_kE,SO(2k-1)]^\bullet\otimes\Q=\pi_{4k-5}(SO(2k-1))\otimes\Q=\Q$. The degree in $[p_kE,SO(2k-1)]^\bullet\otimes \Q$ is defined by the degree in terms of an order 0 generator in $\pi_{4k-5}(SO(2k-1))\cong\Z\oplus\mbox{(finite)}$. (An analogue of Proposition~\ref{prop:[E,SO(5)]}.)

\item $\pi_{4k-5}{SU(4k-4)}/{(SU(2k-1)\times SU(2k-3))}=\Z_{(2k-3)!}$. Hence $\pi_{4k-5}SU(2k-1)\to\pi_{4k-5}SU(4k-4)/SU(2k-3)$ sends the generator to $(2k-3)!$ times the generator. (An analogue of Lemma~\ref{lem:SU->V}.)

\item $p_{k-1}(p_kE;\psi(G)\circ\tau_{E^\bullet})-p_{k-1}(p_kE;\tau_{E^\bullet})=a_{k-1}(2k-3)!\,\deg{G}$ where $a_n=1$ if $n\equiv 0\ \mbox{(mod 2)}$ and $a_n=2$ if $n\equiv 1\ \mbox{(mod 2)}$. Here a result in \cite{Lun} is used. (An analogue of Lemma~\ref{lem:p'2}.)

\item $\zeta_2(p_kE;\psi(G)\circ\tau_{E^\bullet})-\zeta_2(p_kE;\tau_{E^\bullet})=\displaystyle\frac{a_{k-1}(2k-3)!}{48}\,\deg{G}[\Theta]$. (An analogue of Lemma~\ref{lem:dzeta=d2} and \ref{lem:d=theta}.)

\item $L_{k-1}(p_1,\ldots,p_{k-1})=\displaystyle\frac{2^{2k-2}(2^{2k-3}-1)B_k}{(2k-2)!}p_{k-1}+\mbox{(terms of $p_1,\ldots,p_{k-2})$}.$

\item $H^{4p}(E^\bullet\times I,\partial(E^\bullet\times I);\Z)\wedge H^{4(k-1-p)}(E^\bullet\times I,\partial(E^\bullet\times I);\Z)=0$ unless $p=0, k-1$.
\end{itemize}
\end{proof}

%\clearpage
%%%%%%%%%%%%%%%%%%%%%%%%%%%%%%
%%%%%%%%%%%%%%%%%%%%%%%%%%%%%%
%%%%%%%%%%%%%%%%%%%%%%%%%%%%%%
\section{Clasper-bundles}\label{ss:clasper-bundles}

In this section, suspension of graph clasper is defined. Then for a 7-dimensional homology sphere $M$, we shall construct many smooth framed $M$-bundles associated to trivalent graphs, what we will call graph clasper-bundles, by using suspensions. We will show that they are in some sense dual to the Kontsevich classes, which implies the non-triviality of the Kontsevich classes.

More precisely, we shall construct a linear map
\[ \psi_{2n}:\calG_{2n}\to H_{4n}(\BDiff{M};\R) \]
by using `higher-dimensional suspended clasper' construction, and will prove the following 
\begin{Thm}\label{thm:B}
Let $k=4$ and let $M$ be a 7-dimensional homology sphere, then 
\begin{enumerate}
\item The diagram
\[ \xymatrix{
	\calG_{2n} \ar[r]^{\kern-1cm\psi_{2n}} \ar[d]_{\mathrm{proj.}} & H_{4n}(\BDiff M;\R) \ar[d]^{\zeta_{2n}}\\
	\calA_{2n} \ar[r]^{\times 2^{2n}} & \calA_{2n} } \]
is commutative. 
\item $\Im\psi_{2n}$ is included in the image of the Hurewicz homomorphism $\pi_{4n}\BDiff{M}\otimes\Q\to H_{4n}(\BDiff M;\R)$.
\end{enumerate}
\end{Thm}

Composed with any linear map $\calA_{2n}\to\R$, the power $\zeta_{2n}^p$ of $\zeta_{2n}$ for any non-negative integer $p$ yields $\R$-valued characteristic classes. Recall that the degree of a trivalent graph is the number of vertices.
\begin{Cor}\label{cor:dimA}
Suppose that $k=4$ and that $M$ is a 7-dimensional homology sphere.
\begin{enumerate}
\item The degree $2n$ part of $\exp{\sum_{j=1}^\infty\zeta_{2j}}$ yields $\dim\R[\calA_2,\calA_4,\ldots,\calA_{2n}]^{(\deg{2n})}$ linearly independent $\R$-valued characteristic classes of degree $4n$ where $\R[\calA_2,\calA_4,\ldots,\calA_{2n}]$ is the polynomial ring generated by elements of $\calA_2,\calA_4,\ldots,\calA_{2n}$.
\item $\dim\Im\psi_{2n}\geq\dim\calA_{2n}$.
\end{enumerate}
\end{Cor}

\begin{Rem}\label{rem:zdeeper}
The dimensions of the spaces $\calA_{2n}$ for degrees up to 22 are computed in \cite{BN} as follows:
\begin{center}
\begin{tabular}{c|cccccccccccc}\hline
degree ($2n$) & 0 & 2 & 4 & 6 & 8 & 10 & 12 & 14 & 16 & 18 & 20 & 22\\\hline
$\dim{\calA_{2n}}$ & 0 & 1 & 1 & 1 & 2 & 2 & 3 & 4 & 5 & 6 & 8 & 9\\
$\dim{\R[\calA_2,\calA_4,\ldots]^{(\deg\,2n)}}$
	& 1 & 1 & 2 & 3 & 6 & 9 & 16 & 25 & 42 & 50 & 90 & 146\\\hline
\end{tabular}
\end{center}

\end{Rem}

\begin{Cor}
For $n\geq 2$, we have
\[ \begin{split}
	\dim\pi_{4n-1}\Diff(D^7\mbox{ rel }\partial)\otimes\Q
	&=\dim\pi_{4n}B\Diff(D^7\mbox{ rel }\partial)\otimes\Q\\
	&\geq\dim\calA_{2n}.
	\end{split} \]
\end{Cor}
\begin{proof}
Corollary~\ref{cor:dimA} implies that $\dim{\pi_{4n}\BDiff{(D^7\mbox{ rel }\partial)}\otimes\Q}\geq \dim\calA_{2n}$. Further, one can show that {\it if $n\geq 2$ and if $\tau$ and $\tau'$ be two different vertical framings on $E$ that coincide on $E_{q_0}^\bullet$, then $[(E,\tau)]=[(E,\tau')]$ in $\pi_{4n}\BDiff{M}\otimes\Q$.} Indeed, by a similar argument as in Lemma~\ref{lem:H(E,dE)}, we can prove that 
\[ H_i(E^\bullet,\partial E^\bullet\cup E^\bullet_{q_0};\Z)\cong\left\{
	\begin{array}{ll}
		\Z & \mbox{if $i=4n+7$}\\
		0 & \mbox{if $0\leq i\leq 4n+6$}
	\end{array}\right. \]
Recall that $E^\bullet$ denotes the bundle obtained from $E$ by restricting the fiber to $M^\bullet$. Thus we have $[E,SO(7)]^\bullet\otimes\Q\cong\pi_{4n+7}(SO(7))\otimes\Q$ where $[E,SO(7)]^\bullet$ denotes the set of homotopy classes of continuous maps $(E^\bullet,\partial E^\bullet\cup E^\bullet_{q_0})\to (SO(7),1)$. We have $[E,SO(7)]^\bullet\otimes\Q=0$ because it is known that $\pi_{4n+7}(SO(7))$ is finite if $n\geq 2$. So there exists a positive integer $p$ such that $p(E,\tau)$ is equivalent to $p(E,\tau')$. Therefore $[(E,\tau)]=[(E,\tau')]$ in $\pi_{4n}\BDiff{M}\otimes\Q$.

Then it follows that
\[ \pi_{4n}\BDiff{M}\otimes\Q=\pi_{4n}B\Diff(M^\bullet\mbox{ rel }\partial)\otimes\Q\]
for $n\geq 2$, thus we also have
\[ \dim\pi_{4n}B\Diff(D^7\mbox{ rel }\partial)\otimes\Q\geq\dim\calA_{2n}\]
for $n\geq 2$.
\end{proof}

Note that when $k=4$ the correction term for $\hat{\zeta}_2$ for each $E$ defines a bijection between the set of homotopy classes of vertical framings and the set $\{5m+\delta(E)\,|\, m\in\Z\}$ for some $0\leq\delta(E)<5$. Then $\delta$ can be considered as a $\Q/5\Z$-valued homotopy invariant of bundles.
\begin{Cor}\label{cor:z2-nontriv}
Suppose that $k=4$ and that $M$ is a 7-dimensional homology sphere. Then either of the following hold:
\begin{itemize}
\item $\hat{\zeta}_2$ (in Theorem~\ref{thm:AA}) is non-trivial or
\item $\delta$ is non-trivial.
\end{itemize}
\end{Cor}
Proof of Corollary~\ref{cor:z2-nontriv} will be given after the proof of Theorem~\ref{thm:B}.
%%%%%%%%%%%%%%%%%%%%%%%%%%%%%%
%%%%%%%%%%%%%%%%%%%%%%%%%%%%%%
\subsection{Suspended claspers}

Now we define some notions which are generalizations of Habiro's clasper defined in \cite{Hab, Hab2} and which will be elementary pieces in the constructions below. For the details about higher dimensional claspers, see \cite{W}, though we will describe here self-contained definitions of them.

\subsubsection{$I$-claspers}
\begin{figure}
\psfrag{p}[cc][cc]{$p$}
\psfrag{q}[cc][cc]{$q$}
\psfrag{lk=1}[cc][cc]{$\mathrm{Lk}(S^p,S^q)=\pm 1$}
\psfrag{Rp+1}[cc][cc]{$\R^{p+1}$}
\psfrag{Rq+1}[cc][cc]{$\R^{q+1}$}
\psfrag{Sp}[cc][cc]{$S^p$}
\psfrag{Sq}[cc][cc]{$S^q$}
\fig{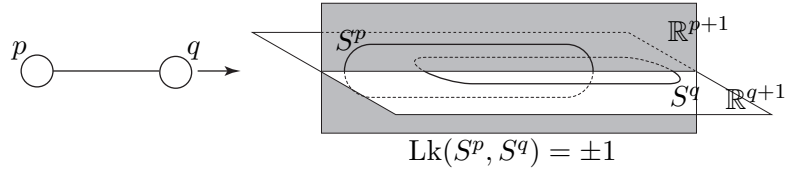}
\caption{$I_{p,q}$-clasper and the associated Hopf link}\label{fig:ipq-link}
\end{figure}

An {\it $I_{p,q}$-clasper} is a normally framed null-homotopic embedding of two disjoint spheres $S^p\sqcup S^q\subset M ^{p+q+1}$ with $p, q\geq 1$ connected by an arc, equipped with a trivialization of the normal $SO(p+q)$-bundle over the arc for which the first $p$-frame is parallel to the $p$-sphere near the one side of the arc and the last $q$-frame is parallel to the $q$-sphere near the other side of the arc. We call each of the two spheres a {\it leaf} and call the arc an {\it edge}. With the given normal framing, we can canonically associate a normally framed two component link to an $I_{p,q}$-clasper by replacing the $I_{p,q}$-clasper with an embedded Hopf link as in Figure~\ref{fig:ipq-link} so that the $p$-sphere lies in the $(p+1)$-plane spanned by the first $p$-frame in the normal frame and by the vector parallel to the direction of the edge, and the $q$-sphere lies in the $(q+1)$-plane spanned by the last $q$-frame in the normal frame and by the vector parallel to the direction of the edge. We orient the two leaves so that the linking number $\mathrm{Lk}(S^p,S^q)$ of the associated Hopf link is 1 if both $p$ and $q$ are odd. By a {\it surgery} along an $I_{p,q}$-clasper, we mean a surgery along its associated framed link.

Since the map induced by the inclusion $\pi_1(SO(2))\to\pi_1(SO(p+q))\cong\Z_2$ or $\Z$, is onto, we can represent the framed edge by an $SO(2)$-framed edge in an untwisted 3-dimensional neighborhood of the edge. This allows us to depict an $I_{p,q}$-clasper in a plane diagram.
One may check that if a $p$-disk and a $q$-disk (or their thickenings) link with the $q$- and $p$-dimensional leaves respectively in a standard way, i.e., each of them intersect the disk bounded by each leaf at a point, then surgery along the $I_{p,q}$-clasper changes partially the linking number of the $p$- and $q$-disk by $\pm 1$. 
%%%%%%%%%%%%%%%%%%%%%%%%%%%%%%
%%%%%%%%%%%%%%%%%%%%%%%%%%%%%%

\subsubsection{Suspended claspers}
Consider a smooth fiber bundle $E\to B$ with fiber a pair $(M, \phi)$ where $\phi$ is a smooth family of $B$-parametrized embeddings of $I_{p,q}$-claspers into $M$ such that it becomes a trivial $M$-bundle if we forget $\phi$. We will call such a bundle a {\it suspended claspers over $B$}.

%%%%%%%%%%%%%%%%%%%%%%%%%%%%%%
%%%%%%%%%%%%%%%%%%%%%%%%%%%%%%

Further we extend the notion of surgery to suspended claspers. Suppose that a suspended clasper $E\to B$ embedded into $M \times B$ forms a trivial sub bundle. Then simultaneous surgery along a suspended clasper, i.e., attaching of $\mbox{(handles)}\times B$ followed by smoothing of corners, yields a possibly non-trivial smooth $M$-bundle. A {\it clasper-bundle} is an $M$-bundle obtained by a sequence of surgeries along suspended claspers.

%\clearpage
%%%%%%%%%%%%%%%%%%%%%%%%%%%%%%
%%%%%%%%%%%%%%%%%%%%%%%%%%%%%%
\subsection{Graph claspers}

Now we review the definition of a higher dimensional graph clasper, which in the case of 3-dimension was also first introduced by Habiro in \cite{Hab}. See \cite{W} for details\footnote{As mentioned in \cite{W}, the definition of the higher dimensional (unsuspended) graph clasper was suggested to the author by Kazuo Habiro, after the author's \cite{W2}.}. Graph clasper itself is not necessary to define graph clasper-bundles below. But it motivates the definition of the graph clasper-bundle. 

\begin{figure}
\psfrag{x}[cc][cc]{$x$}
\psfrag{y}[cc][cc]{$y$}
\psfrag{z}[cc][cc]{$z$}
\fig{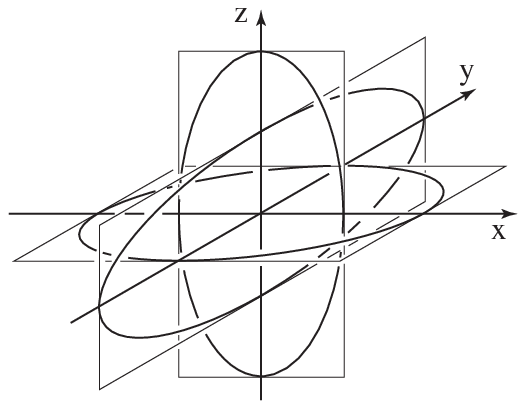}
\caption{}\label{fig:h-borromean}
\end{figure}

In \cite{Hab, Hab2}, the Borromean rings in 3-dimension plays an important role. In the case of higher dimension, the higher dimensional Borromean rings play a similar role. When three natural numbers $p,q,r>0$ satisfy
\begin{equation}\label{eq:borromean-condition}
p+q+r=2m-3,
\end{equation}
one can form higher dimensional Borromean rings $S^p\sqcup S^q\sqcup S^r\to\R^{m}$ as follows. Let $p',q',r'$ be integers such that $p+p'=m-1, q+q'=m-1, r+r'=m-1$. Then $p'+q'+r'=m$. Consider $\R^m$ to be $\R^{p'}\times \R^{q'}\times \R^{r'}$. Then the union of the subsets
\begin{equation}\label{eq:def-borromean}
 \left\{\begin{array}{l}
	\displaystyle S_p\eqdef\{(x,y,z)\in\R^m\,|\,{|x|^2}/{4}+|y|^2=1,\,z=0\}\cong S^p\\
	\displaystyle S_q\eqdef\{(x,y,z)\in\R^m\,|\,{|y|^2}/{4}+|z|^2=1,\,x=0\}\cong S^q\\
	\displaystyle S_r\eqdef\{(x,y,z)\in\R^m\,|\,{|z|^2}/{4}+|x|^2=1,\,y=0\}\cong S^r
	\end{array}\right.
\end{equation}
of $\R^m$ forms a non-trivial 3 component link (see Figure~\ref{fig:h-borromean}). Non-triviality of this link can be proved for example by computing the Massey product of its complement.

Fix an integer $n\geq 3$. A {\it modelled graph clasper} is a connected uni-trivalent graph with
\begin{enumerate}
\item vertex orientation on each trivalent vertex, namely, choices of orders of three incident edges to each trivalent vertex, modulo even number of swappings,
\item decomposition of each edge into a pair of half edges,
\item a natural number $p(h)$ on each half edge $h$ so that if $e=(h_0,h_1)$ is a decomposition of an edge $e$, $p(h_0)+p(h_1)=m-1$ and if $p=p(h_1),q=p(h_2),r=p(h_3)$ are numbers of three incident half edges of a trivalent vertex, then they satisfy the condition (\ref{eq:borromean-condition}),
\item a $p(h_v)$-sphere attached to each univalent vertex $v$ where $h_v$ is the half edge containing $v$.
\end{enumerate}

A {\it graph clasper} is a framed embedding of a modelled graph clasper into an $m$-dimensional manifold together with structures (vertex orientations, $p(\cdot)$). A framed link associated with a graph clasper $G$ is a normally framed link in a regular neighborhood of $G$ obtained by replacing each edge labeled $(p,p')$ with a Hopf link associated to an $I_{p,p'}$-clasper so that the three spheres grouped together at a trivalent vertex form a Borromean rings. Here vertex orientations are used to determine the `orientations' of the Borromean rings (if $m=3$, the Borromean rings $\bar{L}$ obtained from another Borromean rings $L$ by the involution $x\mapsto -x$ in $\R^m$ is not equivalent to $L$).

\begin{figure}
\psfrag{S1}[cc][cc]{$S^1$}
\psfrag{S3}[cc][cc]{$S^3$}
\fig{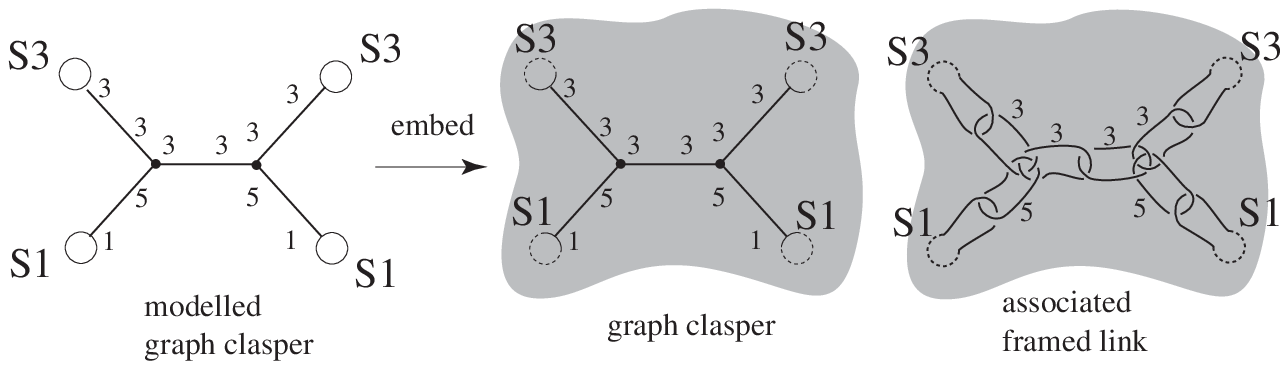}
\caption{}\label{fig:graph-clasper}
\end{figure}

\begin{Exa}
An obvious example is a graph clasper without trivalent vertices. This is just a model of $I_{p,q}$-claspers. Another example of a graph clasper for $m=7$ is depicted in Figure~\ref{fig:graph-clasper}.
\end{Exa}

One may check that graph clasper with cycles exist only if the label $p(h)=1$ is allowed. This condition is always satisfied when $m=3$ or $4$. In the case $m\geq 5$, it may happen that $p(h)>1$ for all $h$. So in that case, graph claspers with cycles do not exist, that is, only the tree shaped graph claspers exist. 

In the case $m=3$, there are many graph claspers so that any trivalent graph yields a graph clasper. However, in the case $m\geq 4$, no trivalent graph yields a graph clasper! In order to construct `dual' objects to the Kontsevich classes for trivalent graphs in high dimensions, we suspend claspers as in the next subsection.

%\clearpage
%%%%%%%%%%%%%%%%%%%%%%%%%%%%%%
%%%%%%%%%%%%%%%%%%%%%%%%%%%%%%
\subsection{Graph clasper-bundles}

We shall define graph clasper-bundles here. More precisely, the goal of this subsection is to define the announced homomorphism $\psi_{2n}$ at the beginning of this section. Let $m=2k-1\geq 3$ be an odd integer. In the following, we restrict only to the $I_{k-1,k-1}$-claspers in $m$-dimensional manifolds for simplicity. 

\subsubsection{Particular suspended three component link}
The following claim is the key observation motivating the definition of the graph clasper-bundle. By an {\it almost $B$-parametrized embedding}, we mean a $B$-parametrized family of smooth maps that are embeddings on $B\setminus\{\mbox{the base point}\}$.
\begin{Obs}\label{obs:param}
There exists an almost $S^{k-2}$-parametrized embedding of a trivial 3 component link into an $m$-ball $B^m(2)$ with radius 2:
\[ \phi_t:S^{k-1}\sqcup S^{k-1}\sqcup S^{k-1}\to B^m(2)\subset \R^m,\ t\in S^{k-2}, \]
considered as distributed in a trivial $B^m(2)$-bundle over $S^{k-2}$, such that the locus of their images, projected into a single $B^m(2)$-fiber, is isotopic to a Borromean rings of dimensions $(k-2,k-2,2k-3)$.
\end{Obs}
For usual graph claspers in \cite{Hab, W} and in the previous subsection, the Borromean rings may be inserted at trivalent vertices. For the definition of the graph clasper-bundles, we will use the `suspended' Borromean rings $\{\phi_t\}_t$ (with a little modification) at trivalent vertices.
\begin{proof}
For an $S^{k-2}$-parametrized 3-component link embedding $\phi_t$, let $\phi_t^{(i)}:S^{k-1}\to B^m(2),\, i=1,2,3$, denote $\phi_t$ restricted to each component. Since the triple $(k-1,k-1,2k-3)$ for $m=2k-1$ satisfies the condition (\ref{eq:borromean-condition}), we can form a Borromean rings $\phi_L$ in $B^m(2)$ of dimensions $(k-1,k-1,2k-3)$ as in the previous subsection. The $(2k-3)$-sphere $L_3$ in $\Im{\phi_L}$ can be considered as a $(k-2)$-fold loop suspension of a $(k-1)$-sphere. Namely, the $(2k-3)$-sphere $L_3$ is covered just once by the locus of an almost $S^{k-2}$-parametrized embedding $\tilde\phi_t$ of $(k-1)$-spheres. Therefore, $\phi_t^{(i)}=\phi_L^{(i)}$ (constant over $t$) for $i=1,2$, and $\phi_t^{(3)}=\tilde\phi_t$ ($t\in S^{k-2}$) gives the desired distribution.
\end{proof}
%\clearpage

\subsubsection{Surgery along the suspended three component link}\label{ss:s-3comp}
Now we want to define correctly a surgery along such a three component parametrized link. In order for such surgery to be well-defined, we are left with the following matters to be overcame:
\begin{enumerate}
\item The image of the almost parametrized embedding $\tilde\phi_t$ defined above degenerates into a point in the fiber of the base point $t^0$ of $S^{k-2}$.
\item We need to prove that the image of $\{\tilde\phi_t\}_t$ after a suitable modification forms a trivial $S^{k-1}$-bundle over $S^{k-2}$ and that it is stabilized near the base point of the fiber $S^{k-1}$.
\end{enumerate}
To overcome these matters, we define a parametrized embedding $\tilde{\varphi}_t:S^{k-1}\to B^m(2)$ by modifying $\tilde{\phi}_t$ so that it is non-degenerate everywhere over $S^{k-2}$. 

\begin{figure}
\psfrag{Qm}[cc][cc]{$Q^m$}
\psfrag{z}[cc][cc]{$S^{k-1}$}
\psfrag{t}[cc][cc]{$S^{k-2}$}
\fig{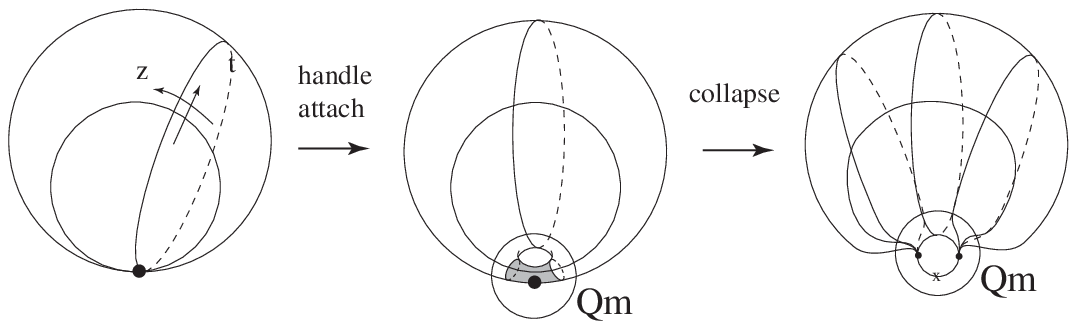}
\caption{}\label{fig:glue-w}
\end{figure}

Let $Q^m\subset B^m(2)$ be an embedded small ball including the base point of $L_3$, where $L_3$ is the third component of the image of $\phi_L$, appeared in Observation~\ref{obs:param}. First we make an embedding of $S^{k-2}\times S^{k-1}$ in $B^m(2)$ by attaching a small $(k-1)$-handle to the $(2k-2)$-disk bounded by $L_3$ along the trivially embedded $(k-2)$-sphere on $L_3$ (see Figure~\ref{fig:glue-w}). Then we collapse the $(k-1)$-handle into its core $(k-1)$-disk so that 
\begin{itemize}
\item the (limiting) boundary of the resulting object is a smooth embedding outside the collapsed part, and
\item the collapsed image from $\{t\}\times S^{k-1}\subset S^{k-2}\times S^{k-1}$ for each $t\in S^{k-2}$ is a smooth embedding.
\end{itemize}
Here we assume that all the changes are included in $Q^m$. Then the resulting family of embeddings by the above construction is the desired one and we will denote it by $\tilde\varphi_t$. See Figure~\ref{fig:glue-w} for an explanation of this construction.

\begin{Prop}\label{prop:Ysus-Y}
The parametrized embedding $(\phi_L^{(1)},\phi_L^{(2)},\tilde\varphi_t)$ can be obtained (up to isotopy) by surgery along a (unsuspended) $Y$-graph clasper in $B^m(2)$ from the trivial one $(\phi_L^{(1)},\phi_L^{(2)},\phi_0^{(3)})$ where 
\begin{itemize}
\item the $Y$-graph clasper is associated with the Borromean rings of dimensions $(2k-3, 2k-3, 2k-3)$,
\item $\phi_0^{(3)}:S^{k-1}\to B^m(2)$ is a constantly parametrized embedding disjoint from $\phi_L^{(1)}$ and $\phi_L^{(2)}$.
\end{itemize}
\end{Prop}
\begin{proof}
After a suitable isotopy, one can push most of $\Im\tilde\varphi_t\subset B^m(2)\times S^{k-2}$ into the fiber $B^m(2)_{t^0}$ of the base point $t^0\in S^{k-2}$. Then the image of $(\phi_L^{(1)},\phi_L^{(2)},\tilde\varphi_t)$ restricts in $B^m(2)_{t^0}$ to a Borromean rings of dimensions $(k-1,k-1,2k-3)$, with something little change near the base point of the third component, that is disjoint from all other components. Then the first two components trivially suspended over $S^{k-2}$ together with the modified $(2k-3)$-sphere in $B^m(2)_{t^0}$, may be seen as a part of the Borromean rings of dimensions $(2k-3,2k-3,2k-3)$ in $(3k-3)$-dimension. Hence the result follows.
\end{proof}
%\clearpage

\subsubsection{Graph clasper-bundle}
We denote by $\phi_t^Y$ the parametrized embedding 
\[ (\phi_L^{(1)},\phi_L^{(2)},\tilde\varphi_t):S^{k-1}\sqcup S^{k-1}\sqcup S^{k-1}\to B^m(2),\quad t\in S^{k-2}\]
defined above. Note that $\phi_t^Y$ can be chosen so that each base point of $S^{k-1}$ is fixed. By using this parametrized embedding, we shall construct graph clasper-bundles.

Let $V$ be an $m$-dimensional handlebody obtained from an $m$-disk by attaching three $(k-1)$-handles along 3 component trivial framed link in the boundary of the $m$-disk. First we shall define the $(V\mbox{ rel }\partial)$-bundle $\pi^Y:V^Y\to S^{k-2}$.

\begin{figure}
\psfrag{M}[cc][cc]{$M$}
\psfrag{G(G)}[cc][cc]{$G(\Gamma)$}
\psfrag{Bm}[cc][cc]{$B^m(2)$}
\psfrag{V1}[cc][cc]{$V_1$}
\psfrag{V2}[cc][cc]{$V_2$}
\fig{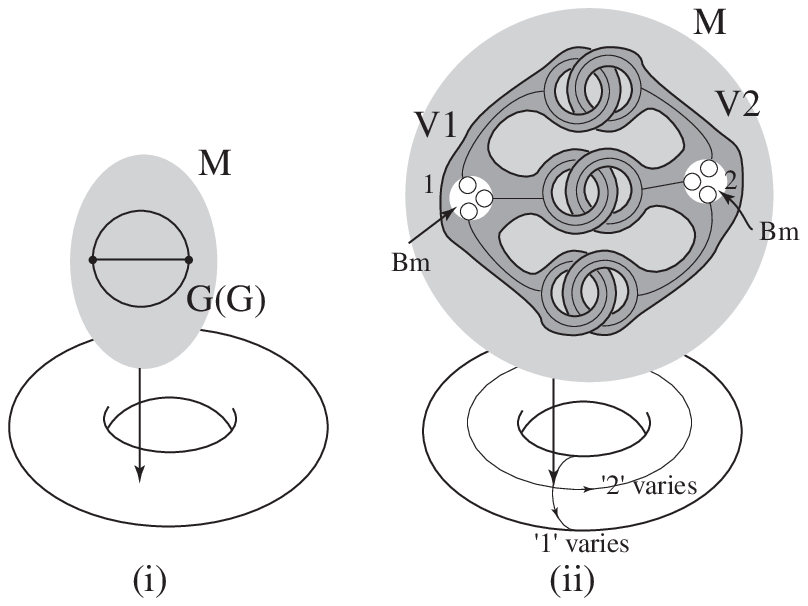}
\caption{}\label{fig:def-clasper-bundle}
\end{figure}

Let us assume that $B^m(2)$ is embedded into the interior of $V$. Then we make a direct product $(V,B^m(2))\times S^{k-2}$ to obtain a trivial sub $B^m(2)$-bundle $\widehat{B}^m(2)\cong B^m(2)\times S^{k-2}$ embedded into the trivial $V$-bundle $V\times S^{k-2}$. Now let 
\[ \phi_I:(I_{k-1,k-1}\sqcup I_{k-1,k-1}\sqcup I_{k-1,k-1})\times S^{k-2}\to V\times S^{k-2}\]
be the three disjoint suspended claspers over $S^{k-2}$ such that 
\begin{enumerate}
\item For each $t\in S^{k-2}$, one of the two leaves of the $i$-th component ($i=1,2,3$) of $\Im{(\phi_I)_t}$ is standardly embedded along the core of the $i$-th $(k-1)$-handle of $V$, and the other leaf is embedded into $V$ isotopically trivial. We denote the latter leaf as an embedding by $(S_i)_t$.

\item $\widehat{B}^m(2)\cap \Im\phi_I$ is precisely a graph (of a function on $t$) of an $S^{k-2}$-parametrized embedding of the three leaves into $B^m(2)$. Thus $\partial \widehat{B}^m(2)\cap \Im\phi_I\cong (\mathrm{pt}\sqcup\mathrm{pt}\sqcup\mathrm{pt})\times S^{k-2}$, the intersection points of edges and leaves of claspers.

\item $(\phi_I)_t$ is standard on $V\setminus B^m(2)$-fiber.

\item $(\phi_I)_t$ restricted to the leaves $S_1\sqcup S_2\sqcup S_3$ coincides with $\phi^Y_t$.
\end{enumerate}
There is an explanation for the form of $\phi_I$ in Figure~\ref{fig:def-clasper-bundle}(ii). Then simultaneous surgeries on $V\times S^{k-2}$ along the suspended claspers $\phi_I$ yield another $(V\mbox{ rel }\partial)$-bundle. We denote the resulting bundle by $\pi^Y:V^Y\to S^{k-2}$.

\begin{Def}[$V^Y$-surgery]
For a given $M$-bundle $\pi:E\to B$, we assume that a trivial $V$-bundle $\widehat{V}\cong V\times B$ is embedded into $E$ as a trivial sub $V$-bundle of $\pi$. Then the {\it $V^Y$-surgery} on $\pi$ along $\widehat{V}$, denoted by $\pi^{V^Y}:E^Y\to B$, is defined for a choice of a $C^\infty$-map $\varphi:B\to S^{k-2}$ as follows:
\[ E^{Y}\eqdef E\setminus\mathrm{int}(\widehat{V})\cup_\partial (\varphi^* V^Y). \]\qed
\end{Def}

\begin{Def}[Graph clasper-bundle]
Let $\Gamma\in\calG$ be a trivalent graph with $2n$ vertices and $3n$ edges not having the part like $\multimap$ and let $G(\Gamma)\subset M$ be a fixed {\it irregular} graph clasper for $\Gamma$ trivially embedded into an $m$-dimensional manifold $M$ with all labels equal to $k-1$. Here `irregular' means that only the condition (\ref{eq:borromean-condition}) for the three labels at trivalent vertices fails to be a graph clasper. Then replace $G(\Gamma)$ with $2n$ disjointly embedded handlebodies $V_1\sqcup \cdots\sqcup V_{2n}$ satisfying the following conditions.

\begin{enumerate}
\item Decompose each handlebody $V_i$ into an $m$-ball $B_i$ and three $(k-1)$-handles $H_1^i\sqcup H_2^i\sqcup H_3^i$ so that $B_i$ includes the $i$-th vertex of $G(\Gamma)$. Then $H_j^i$ and $H_k^{i'}$ are included in $N(e_{ii'})\setminus B_i$, where $j$ and $k$ are determined by the vertex orientation of $G(\Gamma)$, and where $N(e_{ii'})$ denotes a thin tubular neighborhood of the edge of $G(\Gamma)$ connecting the $i$-th and the $i'$-th vertices (if exists). The handles $H_j^i$ and $H_k^{i'}$ link with the linking number $\pm 1$ (1 if $k$ is even), and their positions are determined parallel to the first and the last $p$-frames of the edge.
\item Each edge of $G(\Gamma)$ has just one associated pair of handles $(H_j^i,H_k^{i'})$ as above.
\end{enumerate}
Then the $M$-bundle $\pi^{\Gamma}:E^\Gamma\to (S^{k-2})^{\times 2n}$ is defined as follows. First by taking a direct product $(M, V_1\sqcup \cdots \sqcup V_{2n})\times(S^{k-2})^{\times 2n}$, we obtain $2n$ disjointly embedded trivial sub $V$-bundles in the trivial $M$-bundle $M\times (S^{k-2})^{\times 2n}$. Then we define
\[ E^\Gamma\eqdef \left[ \{M\setminus\mathrm{int}(V_1\sqcup\cdots\sqcup V_{2n})\times B\}\cup_\partial(\varphi_1^*\sqcup\cdots\sqcup \varphi_{2n}^*)V^Y\right] \]
where $\varphi_i:(S^{k-2})^{\times 2n}\to S^{k-2}$ is the $C^\infty$-map for the $V^Y$-sugery along $V_i$, that is the $i$-th projection. We will call such constructed $\pi^{\Gamma}$ a {\it graph clasper-bundle} associated to $\Gamma$. (See Figure~\ref{fig:def-clasper-bundle}.)\qed
\end{Def}

\begin{Rem}
1. The above definition of graph clasper-bundles is also valid for $k=2$, i.e., for graph clasper-bundles consisting of $I_{1,1}$-claspers in a 3-manifold. In this case, the bundle is over $S^0\times\cdots\times S^0$, namely an alternating sum of $Y$-clasper surgeries, which appeared in the context of finite type theory of 3-manifolds \cite{Hab2}.

2. We can generalize the notion of the graph clasper-bundles to arbitrary base $B$ with general choices for $\varphi_i$. In fact there are possibly non isomorphic $V^Y$-surgeries as many as $[B,S^{k-2}]\stackrel{1-1}{\leftrightarrow}\Omega^{\mathrm{fr}}_{\dim{B}-(k-2)}(B)$, the set of bordism classes of normally framed submanifolds of $B$, by the Pontrjagin-Thom construction.
\end{Rem}
%\clearpage

%\begin{figure}
%\psfrag{S}[cc][cc]{$S_v$}
%\psfrag{S'}[cc][cc]{$S'_v$}
%\psfrag{D}[cc][cc]{$D^{k-2}$}
%\fig{stretched-y.eps}
%\caption{}\label{fig:stretched-y}
%\end{figure}

%We will say that a graph clasper-bundle can be {\it quasi vertically framed} if there exists a positive integer $p$ such that $p\pi^{\Gamma}\eqdef\pi^{\Gamma}+\cdots+\pi^{\Gamma}$ ($p$ terms) can be vertically framed.
%\clearpage

\subsubsection{Existence of vertical framings for $k=4$}\label{ss:multiple}
To complete the definition of $\psi_{2n}$, we shall give each graph clasper-bundle a certain vertical framing. 

Let $\pi^\Gamma:E^{\Gamma}\to (S^{k-2})^{\times 2n}$ be a graph clasper-bundle. We make the bundle $\pi^{\Gamma}(2v_i)\ (v_i\in V(\Gamma))$ from $\pi^{\Gamma}$ as follows.

Suppose that the vertex $v_i$ correspond to $V_i$ and let $\varphi_i'\eqdef{\bf 2}\circ\varphi_i:(S^{k-2})^{\times 2n}\to S^{k-2}$ where ${\bf 2}:S^{k-2}\to S^{k-2}$ is the degree 2 map representing twice the generator of $\pi_{k-2}S^{k-2}=\Z$. Then the bundle $\pi^\Gamma(2v_i):E^\Gamma(2v_i)\to (S^{k-2})^{\times 2n}$ is defined similarly as $\pi^\Gamma$ only replacing $\varphi_i$ with $\varphi_i'$. We can apply this construction for several vertices of $\Gamma$ and we will write the result as $\pi^{\Gamma}(2v_{i_1},\ldots,2v_{i_r})$.

The following proposition shows that if $k=4$, any bundle of the form $\pi^{\Gamma}(2v_1,\ldots,2v_{2n})$ is a bundle for which the Kontsevich classes are defined. 
\begin{Prop}\label{prop:EG-framed}
In the case $k=4$, the graph clasper-bundle 
\[ \pi^{\Gamma}(2v_1,\ldots,2v_{2n}):E^{\Gamma}(2v_1,\ldots,2v_{2n})\to(S^{2})^{\times 2n} \]
for any $\Gamma\in\calG_{2n}$ can be vertically framed so that it is standard outside $V_1\sqcup\cdots\sqcup V_{2n}$.
\end{Prop}
The statement given here is stronger than just for saying the existence of the vertical framing because it is used in the proof of Proposition~\ref{prop:induced-form}.
\begin{proof}
Let $E\eqdef E^{\Gamma}(2v_1,\ldots,2v_{2n})$. Assume that $(S^2)^{\times 2n}$ is decomposed into cells obtained from the standard cubic cell decomposition of $(D^2)^{\times 2n}$ by the sequence of collapsings:
\[ D^2\times D^2\times \cdots \times D^2
	\rightarrow S^2\times D^2\times \cdots \times D^2 
	\rightarrow S^2\times S^2\times \cdots \times D^2 
	\rightarrow \cdots \]
where each arrow denotes the collapsing $(D^2,\partial D^2)\to (S^2, \{t^0\})$.

Let $e^2_i$ be the 2-cell corresponding to the $i$-th $S^2$-component of $(S^2)^{\times 2n}$ whose boundary $\partial{e^2_i}$ is to be glued into the base point. We consider the trivial $M$-bundle over $e^2_i$ induced from $\pi^{\Gamma}(2v_1,\ldots,2v_{2n})$ by the inclusion $S^2\hookrightarrow (S^2)^{\times 2n}$. Note that this corresponds to the clasper-bundle for the $Y$-subgraph of $\Gamma$.

We give a polar coordinate on $e^2_i$, namely we identify $e^2_i$ with the set 
\[ \{(r,\theta)\,|\,0\leq r\leq \cos{\theta}, -\frac{\pi}{2}\leq\theta\leq \frac{\pi}{2}\}.\]
For each point $x\in e^2_i$, the diffeomorphism $\varphi_x:E_{q_0}\to E_x$ between the fibers is determined as the result of the smooth deformation along the path $\gamma_x=\{(tr,\theta)\,|\,0\leq t\leq 1\}\ \mbox{for $x=(r,\theta)$}$. Thus we may assume after a homotopy that $\varphi_x=\mathrm{id}$ outside $0\leq\theta\leq \varepsilon$ for some $\varepsilon>0$. Correspondingly, we may assume that the vertical framing is given outside $0\leq\theta\leq \varepsilon$, which is the same as the standard one of $E_{q_0}^\bullet$. On the rest of $e^2_i$, we choose the vertical framing induced by the deformation along path $\gamma_x$. Moreover, since for all $x\in e^2_i$, $\varphi_x$ is identity outside the handlebody $V_i\subset M$, that includes the three $I$-claspers, the vertical framing can be given on $E_x^\bullet\setminus (V_i)_x$ equally to the standard one on $E_{q_0}^\bullet\setminus (V_i)_{q_0}$. 

To show that $\pi^{\Gamma}$ is vertically framed, it suffices to prove the vanishing of the obstructions to homotopy the vertical framing defined above restricted to the trivial $V_i$-bundle $V_i\times\alpha_{\varepsilon}$ over the arc $\alpha_{\varepsilon}=\{(r,\theta)\,|\,r=\cos{\theta}, 0\leq \theta\leq\varepsilon\}$ that is trivialized on $\partial V_i$, into the standard one.  

The obstructions may lie in the following groups:
\[ H^j(V_i\times \alpha_{\varepsilon},\partial (V_i\times \alpha_{\varepsilon});\pi_j(SO(7))),\,0\leq j\leq 8. \]
By Lemma~\ref{lem:homology-VxI} below, we have $H_j(V_i\times\alpha_{\varepsilon},\partial (V_i\times \alpha_{\varepsilon});\Z)=0$ for $0\leq j\leq 4$ and thus the above group is zero for $0\leq j\leq 4$ by the universal coefficient theorem. Further, the above group is zero for $j=5, 6$ because $\pi_5(SO(7))=0, \pi_6(SO(7))=0$. Again by Lemma~\ref{lem:homology-VxI} below, we have $H_j(V_i\times\alpha_{\varepsilon},\partial (V_i\times \alpha_{\varepsilon});\Z)=0$ for $j=6, 7$ and thus the above group is zero for $j=7$. Therefore, the only obstruction may lie in the group
\[ H^8(V_i\times \alpha_{\varepsilon},\partial (V_i\times \alpha_{\varepsilon});\pi_8(SO(7))). \]
Since $\pi_8(SO(7))\cong\Z_2\oplus\Z_2$, the obstruction vanishes after making $\pi^{\Gamma}(2v_i)$.

Since $V_i$'s are mutually disjoint, the vertical framing obtained over the 2-skeleton of $(S^2)^{\times 2n}$ may obviously extends to whole of $(S^2)^{\times 2n}$. Further, since the obtained vertical framing is trivialized on $\partial (V_1\sqcup\cdots\sqcup V_{2n})$, we can extend it to whole of $M$ by the standard vertical framing.
\end{proof}

\begin{Lem}\label{lem:homology-VxI}
Under the settings in the proof of Proposition~\ref{prop:EG-framed}, we have
\[ H_j(V_i\times \alpha_{\varepsilon},\partial(V_i\times \alpha_{\varepsilon});\Z)\cong
	\left\{
	\begin{array}{cl}
		\Z\oplus\Z\oplus\Z & \mbox{if $j=5$}\\
		\Z & \mbox{if $j=8$}\\
		0 & \mbox{otherwise}
	\end{array}
	\right. \]
\end{Lem}
\begin{proof}
By the Poincar\'e-Lefschetz duality, we have
\[ H_j(V_i\times \alpha_{\varepsilon},\partial(V_i\times \alpha_{\varepsilon});\Z)
	\cong H^{8-j}(V_i\times \alpha_{\varepsilon};\Z)
	\cong H^{8-j}(V_i;\Z). \]
\end{proof}

\begin{Rem}
In the proof of Proposition~\ref{prop:EG-framed}, we can choose a cell decomposition of $V_i\times\alpha_\varepsilon$ relative to the boundary, consisting only of three 5-cells and one 8-cell. Since $\pi_5(SO(7))=0$ and $\pi_6(SO(7))=0$, the extension of the null-homotopy to the 7-skeleton is unique up to homotopy. Further by $\pi_9(SO(7))\cong\Z_2\oplus\Z_2$, there are four different extensions of the null-homotopy to the 8-skeleton.
\end{Rem}
%\clearpage

%%%%%%%%%%%%%%%%%%%%%%%%%%%%%%
%%%%%%%%%%%%%%%%%%%%%%%%%%%%%%
\subsection{Duality between graph clasper-bundles and characteristic classes}

Let $k=4$ and let $M$ be a 7-dimensional homology sphere. Let 
\[ \psi_{2n}: \calG_{2n}\to H_{4n}(\BDiff M;\R) \]
be the linear map defined for each connected trivalent graph $\Gamma$ as follows:
\begin{description}
\item[if $\Gamma$ does not have $\multimap$] $\psi_{2n}$ is defined as the class of the image of the classifying map for $\pi^{\Gamma}(2v_1,\ldots,2v_{2n})$ with a choice of a vertical framing $\tau(\Gamma)$ which is standard outside $V_1\sqcup\cdots\sqcup V_{2n}\subset M$. (Such a choice of $\tau(\Gamma)$ is possible by Proposition~\ref{prop:EG-framed}.)
\item[if $\Gamma$ has $\multimap$] $\psi_{2n}$ is defined as $0$.
\end{description}
We will write $[E]$ for the class of the image in $\BDiff{M}$ of the classifying map for a bundle $E\to B$.

%%%%%%%%%%%%%%%%%%%%%%%%%%%%%%
\subsubsection{A choice of the fundamental form in graph clasper-bundles when $k=4$}
The choice of the framing made in Proposition~\ref{prop:EG-framed} allows one to make the fundamental 6-form on $C_2(M)$-bundles more accessible. Let $C(E^{\Gamma})\to B$ be the $C_2(M)$-bundle associated to the $M$-bundle $E^{\Gamma}(2v_1,\ldots,2v_{2n})$ and let $\beta_M\eqdef f^*\alpha_{\Diff{M}}$ where 
\[ f:C(E^\Gamma)\to C_2(M)\gtimes\EDiff{M}\]
is a bundle morphism. To simplify the proof of Theorem~\ref{thm:B}, we replace $\beta_M$ with another one within a cohomology class.

One may check that the value $\langle\zeta_{2n},[E^{\Gamma}(2v_1,\ldots,2v_{2n})]\rangle$ does not change if one replaces the form $\beta_M$ with another form $\beta_M'$ such that 
\begin{itemize}
\item $[(\beta_M)_t]=[(\beta_M')_t]$ in $H^{6}(C_2(M)_t;\R)$ for every $t\in B$,
\item $\iota^*\beta_M'=-\beta_M'$ and 
\item $\beta_M|\partial C_2(M)_t=\beta_M'|\partial C_2(M)_t$ for every $t\in B$.
\end{itemize}
So we shall replace $\beta_M$ with such a $\beta_M'$ so that we can compute the integral explicitly.

For any $i\in\{1,\ldots,2n\}$, fix disjoint simple $S^3$-cycles $(a_j^i)_{j=1,2,3}$ and simple $S^3$-cycles $(b_j^i)_{j=1,2,3}$ on $\partial V_i$ such that
\begin{itemize}
\item $a_j^i$ bounds a 4-disk in $V_i$ and $b_j^i$ bounds a 4-disk in $M\setminus\mathrm{int}(V_i)$.
\item $\langle a_j^i,b_k^i\rangle_{\partial V_i}=\delta_{jk}$.
\end{itemize}
Let $\eta(a_j^i,t)$ be a closed 3-form on $(V_i)_{t=(t_1,\ldots,t_{2n})}$, the $V_i$-fiber over $t$ in $E^\Gamma(2v_1,\ldots,2v_{2n})$, such that its support intersects the thin collar $I\times \partial V_i$ inside $I\times (a_j^i\times D^3)$ where $a_j^i\times D^3$ is a fixed tubular neighborhood of $a_j^i$ in $\partial V_i$, where the restriction of $\eta(a_j^i,t)$ here is the $\varepsilon$-Thom form on $I\times a^j_i$. We can show that $\eta(a_j^i,t)$ can be chosen to be continuous in $t$ (See Appendix~\ref{ss:sim-norm}). The following proposition follows from the observation above and from almost the same proof as \cite[Proposition 3.3]{Les2}. Proof of Proposition~\ref{prop:induced-form} is given in Appendix~\ref{ss:sim-norm}.

\begin{Prop}\label{prop:induced-form}
Suppose that $k=4$ and that the framing of $[E^\Gamma(2v_1,\ldots,2v_{2n})]$ is chosen as $\tau(\Gamma)$ (as in Proposition~\ref{prop:EG-framed}). The form $\beta_M$ on $C(E^\Gamma)$ can be replaced without affecting the resulting value $\langle\zeta_{2n},[E^{\Gamma}(2v_1,\ldots,2v_{2n})]\rangle$ so that:
\begin{itemize}
\item Let $I(t)\subset\{1,\ldots,2n\}$ be the subsets of labels such that $i\in I(t)$ if and only if $t_i\neq t_i^0$. Then For any $t=(t_1,\ldots,t_{2n}), t'=(t_1',\ldots,t_{2n}')\in(S^2)^{\times 2n}$ with $t_i=t_i'\ (\forall i\in I(t)\cap I(t'))$, we have $\beta_M(t_1,\ldots,t_{2n})=\beta_M(t_1',\ldots, t_{2n}')$ on 
\[ C_2\Bigl((M\setminus \bigcup_{i\in I(t)\cup I(t')}\mathrm{int}(V_i))\cup\bigcup_{j\in I(t)\cap I(t')} (V_j)_t\Bigr). \]

\item On $(V_i)_t\times (V_k)_t$,
\[ \beta_M(t)=\sum_{j,k\in\{1,2,3\}}\mathrm{Lk}(b_j^i,b_l^k)\,p_1^*\eta(a_j^i,t)\wedge p_2^*\eta(a_l^k,t). \]
where $p_1, p_2:C_2(M)\to C_1(M)$ denote the first and the second projection, respectively.
\end{itemize}
\end{Prop}
\begin{Rem}\label{rem:eta-t}
In Proposition~\ref{prop:induced-form}, we can choose $\eta(a_j^i,t)$ so that it depends only on $t_i$ (see Appendix~\ref{ss:sim-norm}).
\end{Rem}

\begin{proof}[Proof of Theorem~\ref{thm:B}]
First we assume that the form $\beta_M$ on $C(E^\Gamma)$ has been chosen as in Proposition~\ref{prop:induced-form}.

{\bf (1)} The commutativity of the diagram is a consequence of the following identity:
\[ \langle \zeta_{2n}, [E^{\Gamma}(2v_1,\ldots,2v_{2n})] \rangle = 2^{2n}[\Gamma] \]
for any choice of the vertical framing $\tau({\Gamma})$ that is standard outside $V_1\sqcup\cdots\sqcup V_{2n}$. So we shall now prove this identity.

Let $(t_1,\ldots, t_{2n})$ denote the coordinate of $(S^2)^{\times 2n}$ and let $\omega(\Gamma')(t_1,\ldots,t_{2n})$ be the integrand form for the integral associated to $\Gamma'$, restricted to the configuration space fiber of $(t_1,\ldots, t_{2n})$.

First we see that the computation can be simplified to the one for a bundle with fiber a direct product of some simple spaces. Let $U_i\subset C_{2n}(M)$ be the subset consisting of configurations such that no points are included in $V_i$. We show that the fiber integration restricted to $U_i$-fiber degenerates. We consider the case $i=1$ for simplicity. Let 
\[ \pi_1:S^2\times S^2\times\cdots\times S^2\to \{t_1^0\}\times S^2\times \cdots\times S^2,\quad(t_1^0:\mbox{base point})\]
be the projection defined by $(t_1,t_2,\ldots,t_{2n})\mapsto (t_1^0,t_2,\ldots,t_{2n})$. Then $\pi_1$ can be extended to a bundle morphism $\hat{\pi}_1$ between the sub $U_1$-bundles of $\pi^\Gamma(2v_1,\ldots,2v_{2n})$ and of its restriction to $\{t_1^0\}\times(S^2)^{\times 2n-1}$. Since we can write $\omega(\Gamma')(t_1,t_2,\cdots,t_{2n})=\hat\pi_1^*\omega(\Gamma')(t_1^0,t_2,\cdots,t_{2n})$ over $U_1$ by Proposition~\ref{prop:induced-form}, we have
\[ \begin{split}
	\int_{(t_1,\ldots,t_{2n})\in(S^2)^{\times 2n}}\int_{U_1}\omega(\Gamma')(t_1,t_2,\cdots,t_{2n})
	&=\int_{(S^2)^{\times 2n}}\int_{U_1}\hat\pi_1^*\omega(\Gamma')(t_1^0,t_2,\cdots,t_{2n})\\
	&=\int_{\{t_1^0\}\times (S^2)^{\times 2n-1}}\int_{U_1}\omega(\Gamma')(t_1^0,t_2,\cdots,t_{2n})\\
	&=0
	\end{split} \]
by a dimensional reason. So it suffices to compute the integral over $\widetilde{C}\eqdef C_{2n}(M)\setminus\bigcup_i{U_i}$-fiber. Since at least one point is included in each $V_i$ for any configuration in $\widetilde{C}$, $\widetilde{C}$ is a disjoint union of spaces of the form $V_1\times \cdots \times V_{2n}$.

Further we show that the integration domain can be reduced into a direct product of some bundles as follows. Let $\tilde{V}_i\to S^2$ be the $(V_i\mbox{ rel }\partial)$-bundle induced from $E^{\Gamma}(2v_1,\ldots, 2v_{2n})$ by the inclusion $\iota_i:S^2\hookrightarrow(S^2)^{\times 2n}$ given by $t_i\mapsto (t_1^0,\ldots,t_i,\ldots,t_{2n}^0)$, followed by restriction to the $V_i$-fiber (this is precisely equivalent to $2\pi^Y$). Recall from Remark~\ref{rem:eta-t} that in Proposition~\ref{prop:induced-form}, we can choose $\eta(a_j^i,t)$ so that it depends only on $t_i$. Hence the integral equals
\[ \int_{(S^2)^{\times 2n}}\int_{V_1\times\cdots\times V_{2n}}\omega(\Gamma')
	=\int_{\tilde{V}_1\times\ldots\times \tilde{V}_{2n}}\omega(\Gamma'). \]

In Proposition~\ref{prop:induced-form}, all the $\eta$-forms are standard near $\partial \tilde{V}_i$ and hence the last integral is equal to the integral over $\tilde{V}'_1\times\cdots\times\tilde{V}'_{2n}$, where $\tilde{V}'_i$ denotes the closed manifold obtained from $\tilde{V}_i$ by collapsing $\partial \tilde{V}_i\cong \partial V_i\times S^2$ into $\partial V_i\times\{t_i^0\}$. Thus the integral can be given by a homological evaluation with the fundamental class. Now observe that half of the fundamental homology class of the closed manifold $\tilde{V}'_i$ is represented by the class of the map
\[ \tau_i:H^3(\tilde{V}'_i;\R)\wedge H^3(\tilde{V}'_i;\R)\wedge H^3(\tilde{V}'_i;\R)\to \R \]
corresponding to the triple cup product in $\tilde{V}'_i$ because the suspended $Y$-clasper over an $S^2$ component can be replaced with two disjoint unsuspended $Y$-claspers by Proposition~\ref{prop:Ysus-Y}. Indeed, if $\alpha_i,\beta_i,\gamma_i$ are the classes representing the cores of the three 3-handles of a fiber of $\tilde{V}'_i$ and if $\alpha_i^*,\beta_i^*,\gamma_i^*$ are the duals of $\alpha_i,\beta_i,\gamma_i$ with respect to the evaluation, then
\[ \langle\alpha_i^*\cup\beta_i^*\cup\gamma_i^*,\frac{1}{2}[\tilde{V}'_i]\rangle = 1. \]
Note that $H^9(\tilde{V}'_i;\R)\cong H^3(\tilde{V}'_i;\R)^{\wedge 3}$ is one dimensional and spanned by $\alpha_i^*\wedge\beta_i^*\wedge\gamma_i^*$.

On the other hand, the 6-form $\theta_{e=(i,j)}\eqdef\phi_e^*\beta_M\in\Omega^6(C(E^\Gamma))$ is considered as an element of $H^3(\tilde{V}'_i;\R)\otimes H^3(\tilde{V}'_j;\R)$ corresponding to the linking form. Here $\phi_e$ is defined as in \S\ref{ss:kontsevich-class}.

Therefore, the integral is obtained by contractions of the tensors and we get
\[ \int_{\tilde{V}'_1\times\cdots\times \tilde{V}'_{2n}}
		\tilde{\omega}(\Gamma')
		=\langle \prod_e \theta_{e},
			[\tilde{V}'_1\times\cdots\times \tilde{V}'_{2n}]\rangle
		=\left\{
			\begin{array}{ll}
			|\mathrm{Aut}_e\Gamma|\cdot 2^{2n} & \mbox{if $\Gamma'=\Gamma$}\\
			0 & \mbox{otherwise}
			\end{array}\right.	 
\]
where $|\mathrm{Aut}_e\Gamma|$ is the order of the automorphisms of $\Gamma$ fixing all vertices. See Figure~\ref{fig:contraction} for an explanation of this for the $\Theta$-graph. Here, $\theta_{12}^*=\alpha_1^*\otimes\alpha_2^*+\beta_1^*\otimes\beta_2^*+\gamma_1^*\otimes\gamma_2^*$ and thus 
$\langle {\theta}_{12}^3, [\tilde{V}'_1\times \tilde{V}'_{2}]\rangle
		=3!\cdot\tau_1(\alpha_1^*,\beta_1^*,\gamma_1^*)\cdot\tau_2(\alpha_2^*,\beta_2^*,\gamma_2^*)=3!\cdot 2^2$.
\begin{figure}
\psfrag{t1}[cc][cc]{$\tau_1$}
\psfrag{t2}[cc][cc]{$\tau_2$}
\psfrag{theta}[cc][cc]{$\theta_{12}^*$}
\psfrag{R}[cc][cc]{$\R$}
\psfrag{H3xH3}[ll][ll]{$(H^3(\tilde{V}'_1)\otimes H^3(\tilde{V}'_2))^{\wedge 3}$}
\psfrag{H33xH33}[ll][ll]{$H^3(\tilde{V}'_1)^{\wedge 3}\otimes H^3(\tilde{V}'_2)^{\wedge 3}$}
\fig{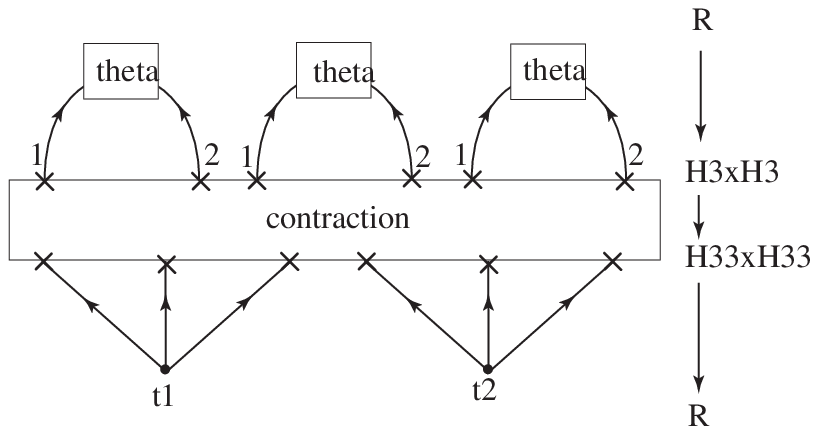}
\caption{}\label{fig:contraction}
\end{figure}

Hence exactly $|\mathrm{Aut}_v\Gamma|\eqdef|\Aut{\Gamma}|/|\mathrm{Aut}_e\Gamma|$ connected components in $\widetilde{C}$ contribute to the term of $\Gamma$ as $2^{2n}$ and the other parts do not contribute. Therefore,
\[ \begin{split}
\zeta_{2n}(E^{\Gamma}(2v_1,\ldots,2v_{2n});\tau(\Gamma))&=|\mathrm{Aut}_v\Gamma|\zeta_{2n}(\tilde{V}'_1\times\cdots\times \tilde{V}'_{2n})\\
	&=|\mathrm{Aut}_v\Gamma|\sum_{\Gamma'}\frac{[\Gamma']}{|\Aut \Gamma'|}\int_{\tilde{V}'_1\times\cdots\times \tilde{V}'_{2n}}
		\tilde{\omega}(\Gamma')\\
	&=|\mathrm{Aut}_v\Gamma|\frac{|\mathrm{Aut}_e\Gamma|\cdot 2^{2n}[\Gamma]}{|\Aut \Gamma|}=2^{2n}[\Gamma]. 
\end{split}\]

{\bf (2)} Observe that $(S^2)^{\times 2n}$ can be made into the one homotopy equivalent to $S^{4n}$ by attaching $2n$ 3-cells along each $S^2$-factor and that the unframed $M^\bullet$-bundle structure extends over the resulting complex $X\simeq S^{4n}$. So we need to consider the obstruction to extend the vertical framing on $E^\Gamma(2v_1,\ldots,2v_{2n})$  over $X$. To do this, we consider the standardly vertical framed trivial $M^\bullet$-bundle $E_i^{\mathrm{cell}}\eqdef M^\bullet\times D^3$ over a 3-cell. Here we may assume that the vertical framing restricted to the boundary of $\tilde{E}_i\eqdef\iota_i^*E^{\Gamma}(2v_1,\ldots,2v_{2n})$ ($\iota_i$: defined above), which is now assumed to be a bundle over the attaching 2-sphere, coincide with that of $E_i^\mathrm{cell}$ restricted to the boundary $\partial D^3$. Moreover, it is not difficult to see that $\tilde{E}_i$ is a trivial $M^\bullet$-bundle as an unframed bundle. We consider the obstruction for the existence of the homotopy between the vertical framings of $\tilde{E}_i$ and of $E_i^\mathrm{cell}|\partial D^3$.

By the Poincar\'{e}-Lefschetz duality, we can show that 
\[ H_j(M^\bullet\times D^2,\partial(M^\bullet\times D^2);\Z)\cong
	\left\{\begin{array}{ll}
		\Z & \mbox{if $j=9$}\\
		0 & \mbox{otherwise}
	\end{array}\right. \]
where $D^2\subset \partial D^3$ is an embedded disk where the obstruction may be included. By $\pi_9(SO(7))\cong\Z_2\oplus\Z_2$, the only obstruction may lie in $H^9(M^\bullet\times D^2,\partial (M^\bullet\times D^2);\pi_9(SO(7)))\cong\Z_2\oplus\Z_2$. So if we replace $E^\Gamma(2v_1,\ldots,2v_i,\ldots,2v_{2n})$ with $E^{\Gamma}(2v_1,\ldots,2^2v_i,\ldots,2v_{2n})$, the vertical framing extends to $E^{\Gamma}(2v_1,\ldots,2^2v_i,\ldots,2v_{2n})\cup_{E_i^\mathrm{cell}|\partial D^3}(M^\bullet\times D^3)$. Therefore, the vertical framing on $E^{\Gamma}(2^2v_1,\ldots,2^2v_{2n})$ can be extended over $X$. Finally, by collapsing the attached 3-cells into the base point by a homotopy, we obtain a vertically framed bundle associated to a class in $\pi_{4n}\BDiff{M}$. Since $[E^{\Gamma}(2^2v_1,\ldots,2^2v_{2n})]=2^{2n}[E^{\Gamma}(2v_1,\ldots,2v_{2n})]$ in $H_{4n}(\BDiff{M};\R)$ and the attaching of 3-cells corresponds to a homotopy in $\BDiff{M}$, the result follows.
\end{proof}

\begin{proof}[Proof of Corollary~\ref{cor:z2-nontriv}]
From the proof of Theorem~\ref{thm:B}, we have $\langle \zeta_2, [E^\Theta(2^2v_1,2^2v_2)] \rangle = 2^2\cdot 2^2[\Theta]$ and that  $[E^\Theta(2^2v_1,2^2v_2)]$ is homologous in $H_{4n}(\BDiff{M};\R)$ to an image from an element of $\pi_{4n}\BDiff{M}$. Recall from the proof of Theorem~\ref{thm:AA} that the framing dependence of $\zeta_2$ for a change of framing $G\in [E^\Theta(2^2v_1,2^2v_2),SO(7)]^\bullet$ is $\displaystyle\frac{2\cdot 5!}{48}\,\deg{G}[\Theta]=5\,\deg{G}[\Theta]$. 

If $\delta=0$, then $\hat{\zeta}_2(E^\Theta(2^2v_1,2^2v_2))$ must be of the form $(16+5\,\deg{G})[\Theta]$. Since $\deg{G}$ is an integer, it follows that $\hat{\zeta}_2(E^\Theta(2^2v_1,2^2v_2))$ is non-zero.
\end{proof}

%%%%%%%%%%%%%%%%%%%%%%%%%%%%%%
%%%%%%%%%%%%%%%%%%%%%%%%%%%%%%
\section{Further directions}\label{s:directions}

Now we shall briefly remark some direction expected to be studied after the present paper. 

In the case of 3-dimensional homology sphere, there is a very powerful theory producing a lot of topological invariants, a theory of finite type invariants, initiated by Ohtsuki in \cite{Oh}. It is conjectured that any different prime homology 3-spheres are distinguished by finite type invariants. Using the construction of Le-Murakami-Ohtsuki of a universal invariant \cite{LMO}, Le proved that the Le-Murakami-Ohtuski invariant is universal among $\R$-valued finite type invariants of homology 3-spheres \cite{Le} and it turned out that there are $\dim\R[\calA_2,\calA_4,\ldots,\calA_{2n}]^{(\deg{\leq 2n})}$ linearly independent $\R$-valued finite type invariants of degree $\leq n$. (Overview of the results related to Ohtsuki's finite type invariant is explained in detail in \cite{Oh2}.)

The construction of the invariant by Le-Murakami-Ohtsuki is based on the Kirby calculus \cite{Kir}. Namely, they use the representation of a 3-manifold by a framed link in $S^3$ considered modulo some moves on them, called the Kirby moves:
\begin{equation}\label{eq:surg}
\xymatrix{
	& \{\mbox{framed links in $S^3$}\}/\sim \ar[r]^{\kern-1cm\mathrm{surgery}} \ar[d]  & \{\mbox{closed ori. connected 3-manifolds} \}/\sim \\
	& \{\mbox{framed links in $S^3$}\}/\mbox{($\sim$, Kirby moves)} \ar[ru]_{\sim} & \\
	}
\end{equation}

 They consider the Kontsevich integral of the framed link \cite{Kon2}, whose value is an infinite linear sum of certain graphs. Then they invented an operator on the space of the graphs so that the resulting value is in $\R[\calA_2,\calA_3,\ldots]$ and they proved that it is invariant under the Kirby moves, namely, that it is a topological invariant. The Le-Murakami-Ohtsuki construction allows an algorithmic construction of 3-manifold invariants because the Kontsevich framed link invariant can be constructed algorithmically (e.g., \cite{BN2}).

Their construction may be explained in other words as follows. It is obvious that any 3-manifold invariants are pulled back by surgery correspondence in (\ref{eq:surg}) to give framed link invariants:
\begin{equation}\label{eq:H0(surg)}
 H^0(B\Diff M;\R)\mapright{\mathrm{surgery}^*} H^0(\mathrm{Emb}_f(S^1\sqcup\ldots\sqcup S^1, S^3);\R)
\end{equation}
where $\mathrm{Emb}_f(A,B)$ denotes the space of normally framed embeddings $A\hookrightarrow B$. Le-Murakami-Ohtsuki's construction is in some sense an inverse of this.

We expect that there is also an algorithmic construction for the Kontsevich classes restricted to some surgery defined bundles. Namely, the higher dimensional analogue of (\ref{eq:H0(surg)}) may be
\begin{equation}\label{eq:Hp(surg)}
 \mbox{``}H^p(B\Diff M;\R)\mapright{\mathrm{surgery}^*} H^p(\widetilde{\mathrm{Emb}}_f^0(S^{p_1}\sqcup\ldots\sqcup S^{p_r}, M);\R)\mbox{"}
\end{equation}
Here the embeddings have to be restricted to the class such that surgery along which do not change the diffeomorphism type of $M$. Then it is natural to expect that some universal characteristic classes of $M$-bundles are obtained from some cohomology classes of the space of link embeddings. Note that for the map (\ref{eq:Hp(surg)}) to be well-defined, the surgery has to be defined in a canonical way. If the space $\widetilde{\mathrm{Emb}}_f^0(S^{p_1}\sqcup\ldots\sqcup S^{p_r},M)$ is replaced with the space of bounding links such that a disk bounding each knot component does not have self-intersection, then the map surgery$^*$ is well-defined. Anyway, we have a well-defined map
\[ \Hom(\Omega_p(B\Diff M),\R)\mapright{\mathrm{surgery}^*} \Hom(\Omega_p(\widetilde{\mathrm{Emb}}_f^0(S^{p_1}\sqcup\ldots\sqcup S^{p_r}, M)),\R). \]

The following problems are related to the algorithmic construction of the Kontsevich classes.
\begin{Prob}
Find an algorithmic construction for cocycles on the link embedding space.
\end{Prob}
As in 3-dimension, algorithmic construction seems easier for the link embedding space than for $B\Diff{M}$.
\begin{Prob}
Give a smooth bundle analogue of the Kirby calculus.
\end{Prob}

By the way, the following problem might be related to the estimation of the cohomology of $B\Diff{M}$.
\begin{Prob}
How general is the class of (bordism classes of) bundles which are obtained by clasper-bundle surgery?
\end{Prob}
%%%%%%%%%%%%%%%%%%%%%%%%%%%%%%
%%%%%%%%%%%%%%%%%%%%%%%%%%%%%%
\begin{appendix}
\section{The closed form $\alpha_{\Diff{M}}$}\label{s:alpha-form}

The construction of the Kontsevich classes requires a `fundamental' closed form $\alpha_{\Diff{M}}$ on $C_2(M)\gtimes\EDiff{M}$. We shall give a proof that there exists such a well-defined closed form $\alpha_{\Diff{M}}$, which is omitted in \cite{Kon}.

From the Serre spectral sequence of the fibration
\[ (C_2(M),\partial C_2(M))
	\to (C_2(M)\gtimes\EDiff{M},
		\partial C_2(M)\gtimes\EDiff{M})
	\to \BDiff{M}, \]
we have the following
\begin{Lem}\label{lem:serre-ss}
There exists a spectral sequence with
\[ \begin{split}
	&E^{p,q}_2\cong H^p(\BDiff{M};\{H^q(C_2(M)_b,\partial C_2(M)_b;\R)\}_{b\in\BDiff{M}})\\
	&\Rrightarrow H^{p+q}(
		C_2(M)\gtimes\EDiff{M},
		\partial C_2(M)\gtimes\EDiff{M}; \R).
	\end{split}\]
\end{Lem}

The following lemma can be proved by exactly the same way as \cite[Lemma 2.1]{Les}.
\begin{Lem}\label{lem:H(C)=H(S)}
$H_*(C_2(M);\Z)\cong H_*(S^{m-1};\Z). $
\end{Lem}

\begin{Lem}\label{lem:H(C,dC)=0}
For any $b\in\BDiff{M}$ and for $0\leq q\leq m$, $H^q(C_2(M)_b, \partial C_2(M)_b;\R)\cong 0$.
\end{Lem}
\begin{proof}
In this proof, all the homology coefficients are assumed in $\R$. By the Poincar\'{e}-Lefschetz duality and Lemma~\ref{lem:H(C)=H(S)}, we have
\[ H^q(C_2(M)_b,\partial C_2(M)_b)\cong H_{2m-q}(C_2(M))\cong H_{2m-q}(S^{m-1})\cong 0\ (0\leq q\leq m). \]
\end{proof}

\begin{Lem}\label{lem:H(C*E,dC*E)=0}
For  $0\leq q\leq m$, $H^q(C_2(M)\gtimes\EDiff{M},
		\partial C_2(M)\gtimes\EDiff{M}; \R)\cong 0$.
\end{Lem}
\begin{proof}
This follows immediately from Lemma~\ref{lem:serre-ss} and Lemma~\ref{lem:H(C,dC)=0}.
\end{proof}

\begin{Lem}\label{lem:H(CE)=H(dCE)}
The inclusion induces an isomorphism
\[ H^{m-1}(C_2(M)\gtimes\EDiff{M};\R)
	\cong H^{m-1}(\partial C_2(M)\gtimes\EDiff{M}; \R).\]
\end{Lem}
\begin{proof}
This follows from the cohomology exact sequence of the pair 
\[(C_2(M)\gtimes\EDiff{M},
		\partial C_2(M)\gtimes\EDiff{M})\]
and from Lemma~\ref{lem:H(C*E,dC*E)=0}.
\end{proof}

Since we can define a closed $(m-1)$-form on $\partial C_2(M)\gtimes\EDiff{M}$ uniquely determined by the framing, there exists a well-defined closed $(m-1)$-form $\alpha_{\Diff{M}}$ on $C_2(M)\gtimes\EDiff{M}$ by Lemma~\ref{lem:H(CE)=H(dCE)}. Note that the vertical framing on $M\gtimes\EDiff{M}$ determines a trivial $S^{m-1}$-bundle structure on $\partial C_2(M)\gtimes\EDiff{M}$ and thus the closed $(m-1)$-form on $\partial C_2(M)\gtimes\EDiff{M}$ is non-trivial in cohomology.
	
\section{Pushforward}\label{s:pushforward}
Let $\pi:E\to B$ be a bundle with $m$-dimensional fiber $F$. Then the {\it push-forward} (or {\it integral along the fiber}) $\pi_*\omega$ of an $(m+p)$-form $\omega$ on $E$ is a $p$-form on $B$ defined by
\[ \int_c \pi_*\omega=\int_{\pi^{-1}(c)}\omega, \]
where $c$ is a $p$-dimensional chain in $B$. 

Let $\pi^{\partial}:\partial_F E\to B$ be the restriction of $\pi$ to $\partial F$-bundle with the orientation induced from $\mathrm{Int}{(F)}$, i.e., $O_{\partial F}=i(n)O_F$ where $n$ is the in-going normal vector field over $\partial F$. Then the generalized Stokes theorem for the pushforward is 
\begin{equation}\label{eq:stokes}
 d\pi_*\omega=\pi_*d\omega+(-1)^{\deg{\pi_*^{\partial}\omega}}\pi^{\partial}_*\omega.
\end{equation}

%%%%%%%%%%%%%%%%%%%%%%%%%%%%%%
%%%%%%%%%%%%%%%%%%%%%%%%%%%%%%
%%%%%%%%%%%%%%%%%%%%%%%%%%%%%%
\section{Simultaneous normalization of the $\beta_M$-forms}\label{ss:sim-norm}

Here we give a simultaneous normalization of the form $\beta_M$ on $C_2(M)$ for 7-dimensional homology spheres, based on the line of a part of \cite[Proposition~3.3]{Les2}. In this section, we denote the fiber over the base point $t^0\in (S^2)^{\times 2n}$ of $E^\Gamma(2v_1,\ldots,2v_{2n})$ by $M$. We identify a regular neighborhood of $\partial V_i\subset M$ with $[-4,4]\times \partial V_i$ and for $s\in [-4,4]$, set 
\[ V_i[s]\eqdef\left\{
	\begin{array}{ll}
		V_i \cup ([0,s]\times \partial V_i) & \mbox{if }s\geq 0\\
		V_i\setminus ((s,0]\times \partial V_i) & \mbox{if }s\leq 0
	\end{array}\right. \]
Let $S(a_j^i)\subset V_i[4]$ and $S(b_k^i)\subset M\setminus \inn(V_i)$ be the 4-disks bounded by $4\times a_j^i$ and $b_k^i$ respectively, such that if $\mathrm{Lk}(a_j^i,a_{j'}^{i'})=1$ for $i\neq i'$, then $S(b_j^i)\cap V_{i'}=S(a_{j'}^{i'})$, and if $\mathrm{Lk}(a_j^i,a_{j'}^{i'})=0$, then $S(b_j^i)\cap V_{i'}=\emptyset$. 

Let $\eta(b_j^i)$ be the closed 3-form supported in an $\varepsilon$-tubular neighborhood $N_\varepsilon S(b_j^i)$ of $S(b_j^i)$ which is restricted to the Thom class in $H^3(N_\varepsilon S(b_j^i)_x,\partial (N_\varepsilon S(b_j^i)_x);\R), x\in S(b_j^i)$, and $\eta(a_j^i)$ is defined by the pullback by the inclusion $N_\varepsilon S(a_j^i)\to N_\varepsilon S(b_{j'}^{i'})$ for some $i', j'$. 

Fix a base point $p^i$ on $\partial V_i$ and let $\omega(p^i)$ be a closed 6-form supported in a tubular neighborhood of the union of the path $[p^i,\infty]$ and $\partial C_1(M)$ such that it restricts as the usual volume form on $\partial C_1(M)=S^6$ and such that the support is disjoint from all $V_i[4]$ and from all the supports of the above forms. First we shall normalize $\beta_M$ on the subset $V_i\times (C_1(M)\setminus V_i[3])\subset C_2(M)$.

\begin{Prop}\label{prop:normalize-part}
For any subset $N\subset \{1,\ldots, 2n\}$, we can choose $\beta_M$ on $C_2(M)$ so that:
\begin{enumerate}
\item For every $i\in N$, the restriction of $\beta_M$ to $V_i\times (C_1(M)\setminus V_i[3])\subset C_2(M)$ equals
\[ \sum_{(j,k)\in\{1,2,3\}} \mathrm{Lk}(b_j^i, a_k^i[4])\,p_1^*\eta(a_j^i)\wedge p_2^*\eta(b_k^i)+p_2^*\omega(p^i) \]
where $p_1, p_2:C_2(M)\to C_1(M)$ denote the first and the second projection, respectively.
\item $\beta_M$ is antisymmetric with respect to $\iota$ and fundamental, that is closed and $\beta_M|\partial C_2(M)=p_M^*\omega_{S^{m-1}}$.
\end{enumerate}
\end{Prop}

Assume Proposition~\ref{prop:normalize-part} for the moment.  Let $E^\Gamma(i)$ be the pullback bundle from $E^\Gamma(2v_1,\ldots,2v_{2n})$ by the inclusion $S^2\hookrightarrow (S^2)^{\times 2n}$ and let $\tilde{V}_i[s]$ be the sub $(V_i[s]\mbox{ rel }\partial)$-bundle of $E^\Gamma(i)$. We extend $\eta(a_j^i)$ and $\eta(b_k^i)$ to the globally defined forms $\eta(a_j^i,t)$ and $\eta(b_k^i,t)$ on $\tilde{V}_i[4]$ and $E^\Gamma(i)\setminus \mathrm{int}(\tilde{V}_i)$ respectively, as follows. 

Observe that there exists a $(4+2)$-manifold $\widetilde{S}(a_j^i)$ included in $\tilde{V}_i[4]$, bounded by $(4\times a_j^i)\times S^2\subset \tilde{V}_i[4]$, such that it restricts to $S(a_j^i)$ in the fiber over $t^0$. Indeed, the third component of the locus of the parametrized link of Observation~\ref{obs:param} bounds a 6-disk if we ignore the other two components. This bounded 6-disk can be considered as a collection of bounded 4-disks parametrized by $t\in S^2$. So this collection can be suspended over $S^2$ with some intersections with the other components. Those intersections can be removed by suitable attachings of handles parallel to the other two components. The resulting 6-manifold is as desired. Then $\eta(a_j^i,t)$ is defined as the restriction of the $\varepsilon$-Thom form over $\widetilde{S}(a_j^i)$ to the fiber of $t$ . $\eta(b_k^i,t)$ may be naturally extended from $\eta(a_j^i,t)$'s by using $\eta(b_k^i)$'s.

For $I\subset\{1,\ldots,2n\}$ and for $t\in (S^2)^{\times 2n}$ such that $I(t)\subset I$, define $\beta_{M_t}^0$ on 
\[ D_I(\beta_{M_t}^0)\eqdef \bigl(C_2(M_t)\setminus \bigcup_{i\in I}(V_i[-1]_t\times V_i[3]_t)\cup(V_i[3]_t\times V_i[-1]_t)\bigr)
	\cup p_{12}^{-1}\Delta_{M_t\setminus\{\infty\}} \]
where $p_{12}:C_2(M_t)\to M_t\times M_t$ be the projection, so that
\begin{itemize}
\item $\beta_{M_t}^0=\beta_M$ on $C_2(M_t\setminus\cup_{i\in I}V_i[-1]_t)=C_2(M\setminus\cup_{i\in I}V_i[-1])$,
\item \[\beta_{M_t}^0=\sum_{(j,k)\in\{1,2,3\}^2}
		\mathrm{Lk}(b_j^i,a_k^i[4])\,p_1^*\eta(a_j^i,t)\wedge p_2^*\eta(b_k^i,t)+p_2^*\omega(p^i)\]
on $p_{12}^{-1}((V_i)_t\times (M_t\setminus V_i[3]_t))$ when $i\in I$.
\item $\beta_{M_t}^0=-\iota^*\beta_{M_t}^0$ on $p_{12}^{-1}((M_t\setminus V_i[3]_t)\times (V_i)_t)$ when $i\in I$.
\item $\beta_{M_t}^0=p_{M_t}^*\omega_{S^{m-1}}$ on $\partial C_2(M_t)$.
\end{itemize}
\begin{figure}
\psfrag{VI}[cc][cc]{$V_i$}
\fig{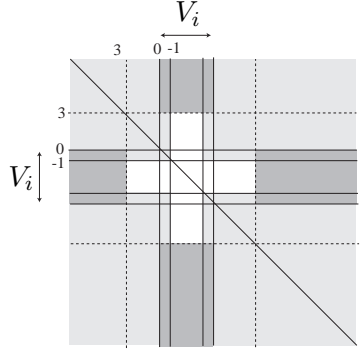}
\caption{Area of $C_2(M)$ where $\beta_M$ is normalized}\label{normalize-area}
\end{figure}
Note that this condition is consistent. In particular, by Proposition~\ref{prop:EG-framed}, the first and the fourth conditions are compatible. Let $C(E^\Gamma(i))\eqdef\cup_{t=(t_1^0,\ldots,t_i,\ldots,t_{2n}^0)}C_2(M_t)$ denote the $C_2(M)$-bundle over $S^2$ associated with $E^\Gamma(i)$. In the following we shall see that the form $\beta_{M_t}^0$ defined over the bundle $D(E^\Gamma(i))\eqdef\cup_{t=(t_1^0,\ldots,t_i,\ldots,t_{2n}^0)}D_{\{i\}}(\beta_{M_t}^0)$ extends to a fundamental 6-form on $C(E^\Gamma(i))$.

\begin{Lem}\label{lem:ssd}
There exists a homology spectral sequence with
\[ E_{p,q}^2\cong H_p(S^2;H_q(D_{\{i\}}(\beta_{M}^0);\R))
	\Rrightarrow H_{p+q}(D(E^\Gamma(i)); \R)
\]
such that $E^2_{p,q}=E^\infty_{p,q}$ if $p+q\leq 6$ and such that $E^\infty_{p,q}= 0$ if moreover $p\notin\{0,2\}$ or $q\notin\{0,4,6\}$. In particular, $H_6(D(E^\Gamma(i));\R)=E^\infty_{0,6}\oplus E^\infty_{2,4}=E^2_{0,6}\oplus E^2_{2,4}$.
\end{Lem}

\begin{Lem}\label{lem:E24}
\begin{enumerate}
\item $E_{2,4}^\infty$ coincides with the kernel of the map induced by the inclusion
\[ H_6(D(E^{\Gamma}(i));\R)\to H_6(C(E^\Gamma(i));\R). \]
\item $\beta_{M_t}^0$ evaluated on $E_{2,4}^\infty$ vanishes.
\end{enumerate}
\end{Lem}

Proofs of Lemma~\ref{lem:ssd} and \ref{lem:E24} will be given later. It follows from these lemmas that the form $\beta^0(i)_t\eqdef \beta^0_{M_t}\,(t=(t_1^0,\ldots,t_i,\ldots,t_{2n}^0))$ on $D(E^\Gamma(i))$ is in the image of the map
\[ H^6(C(E^\Gamma(i));\R)\to H^6(D(E^\Gamma(i));\R). \]
Namely, $\beta^0(i)$ extends to a closed form $\beta^1(i)$ on $C(E^\Gamma(i))$ by the de Rham theorem, and
\[ \beta(i)\eqdef\frac{\beta^1(i)-\iota^*\beta^1(i)}{2} \]
is a fundamental form.

For any $t\in (S^2)^{\times 2n}$, we define
\[ \beta_{M_t}=\left\{
	\begin{array}{ll}
		\beta_{M_t}^0 & \mbox{on }C_2(M_t)\setminus \bigcup_{i\in I(t)}(V_i[-1]_t\times V_i[3]_t)\cup(V_i[3]_t\times V_i[-1]_t)\\
		\beta(i)_t & \mbox{on }C_2(V_i[4]_t)\mbox{ for }i\in I(t)
	\end{array} \right. \]
Then $\beta_{M_t}$ is the required form of Proposition~\ref{prop:induced-form}.

\begin{proof}[Proof of Proposition~\ref{prop:normalize-part}]
We first prove the proposition for the case $N=\{1\}$. Let $\beta_0$ be a fundamental 6-form on $C_2(M)$ and let $\beta$ be the closed 6-form on $V_1[1]\times (C_1(M)\setminus \inn V_1[2])$ defined by the statement. Since integrals for both $\beta_0$ and $\beta$ coincide on $H_6(V_1[1]\times (C_1(M)\setminus \inn V_1[2]);\R)$, there exists a 5-form $\eta$ on $V_1[1]\times (C_1(M)\setminus \inn V_1[2])$ such that 
\[ \beta=\beta_0+d\eta. \]
Here we may assume that $\eta=0$ on $V_1[1]\times \partial C_1(M)$ because $\eta$ is closed on $V_1[1]\times \partial C_1(M)$ and hence exact there.

We further modify $\beta$ so to coincide with $\beta_0$ on $\partial C_2(M)$. Let $\chi$ be a smooth function on $C_2(M)$ supported in $V_1[1]\times (C_1(M)\setminus \inn V_1[2])$, and constant equal to 1 on $V_1\times (C_1(M)\setminus V_1[3])$. Then set
\[ \beta_a\eqdef \beta_0+d(\chi\eta). \]
$\beta_a$ is as required on $V_1\times (C_1(M)\setminus V_1[3])$ and coincides with $\beta_0$ on $\partial C_2(M)$ because $d(\chi\eta)=0$ there. 

Similar modification to $\beta_a$ for $(C_1(M)\setminus V_1[3])\times V_1$, that can be done disjointly from the previous ones, yields another 6-form $\beta_b$ that is as required on
\[ \partial C_2(M) \cup (V_1\times (C_1(M)\setminus V_1[3]))\cup ((C_1(M)\setminus V_1[3])\times V_1). \]
Thus $\beta_M\eqdef(\beta_b-\iota^*\beta_b)/2$ is the required form for $N=\{1\}$.

Now we prove the proposition for general $N$ by induction on $|N|=i$. Let $\beta_0$ be the 6-form satisfying all the hypotheses for $N=\{1,\ldots,i-1\}$, and let $\beta$ be the 6-form satisfying the hypotheses on $\{i\}$ obtained by the first step from $\beta_0$, replacing $V_i$ with $V_i[1]$. Then there exists a 5-form $\eta$ such that $\beta=\beta_0+d\eta$ where $\eta$ may be assumed to vanish on $\partial C_2(M)$ because $H^5(\partial C_2(M);\R)=0$.

Let $\chi$ be a smooth function $\chi$ supported in $V_i[1]\times (C_1(M)\setminus \inn V_i[2])$, that is constant equal to 1 on $V_i\times (C_1(M)\setminus V_i[3])$, and let $\beta_a\eqdef \beta_0+d(\chi\eta)$. Then $\beta_a$ is as required on
\[ \partial C_2(M) \cup \bigcup_{k\in N}(V_k\times (C_1(M)\setminus V_k[3]))
				\cup \bigcup_{k\in N\setminus \{i\}}((C_1(M)\setminus V_k[3])\times V_k). \]
So we need to prove that $\beta_a$ is as required in $V_i[1]\times (\partial C_1(M)\cup \bigcup_{k=1}^{i-1}V_k)$, where the support of $\chi$ intersects the previous changes for $\beta_0$. By the assumptions, $\eta$ may be assumed to vanish on $V_i[1]\times \partial C_1(M)$ and is closed on $V_i[1]\times V_k$ for $i\neq k$. Further by $H^5(V_i[1]\times V_k;\R)=0$, we may assume that $\eta$ vanishes on $V_i[1]\times V_k$.

Finally, by similar modifications as in the first step, we can modify $\beta_a$ so that it integrates correctly as required, and antisymmetric with respect to $\iota^*$.
\end{proof}

%\begin{proof}[Proof of Lemma~\ref{lem:normalize-part-2}]
%In each inductive step of the proof of Proposition~\ref{prop:normalize-part}, one may replace $\eta$ with $\eta+\eta_c$ where $\eta_c$ a linear combination of the closed forms $\eta(b_j^i,t)$, without changing the required properties. Since this replacement adds 
%\[ \int_{p^i\widetilde{\times}([2,3]\times a_j^i\times S^2)}d(\chi\eta_c)=\int_{p^i\widetilde{\times}(3\times a_j^i\times S^2)}\eta_c \]
%to $\int_{p^i\widetilde{\times}\widetilde{S}(a_j^i)}\tilde\beta_a$, one can make $\int_{p^i\widetilde{\times}\widetilde{S}(a_j^i)}\tilde\beta_a$ vanish by adding a suitable $\eta_c$. Since this integral can be considered as a linear functional on $\span\{a_1^i,a_2^i,a_3^i\}$, one can choose $\eta_c$ so that all the integral of the above form vanish simultaneously.
%\end{proof}

\begin{proof}[Proof of Lemma~\ref{lem:ssd}]
First we compute the homology of $D_{\{i\}}(\beta_M^0)$. For any submanifold $X$ of $M$, we denote by $STX$ the face of $\partial C_2(X)$ corresponding to the blow up along the main diagonal $\Delta_X\subset X^{\times 2}$. Since the inclusion from $D_{\{i\}}(\beta_{M}^0)$ to $(C_2(M)\setminus C_2(V_i[-1]))\cup STV_i$ is a homotopy equivalence, it suffices to compute the homology of the latter space. 

Let $\overline{M}=C_1(M)$ and $V=V_i$. We compute the homology of $C_2(M)\setminus C_2(V)\simeq \breve{C}_2(M)\setminus \breve{C}_2(V)$ where $\breve{C}_2(X)\eqdef X^{\times 2}\setminus\{\mbox{diagonal}\}$. Observe that 
\[ H_*(\cM\setminus V)=\left\{
	\begin{array}{ll}
		\R[\partial \cM] & \mbox{if $*=6$}\\
		\R[a_1^i[4]]\oplus\R[a_2^i[4]]\oplus\R[a_3^i[4]] & \mbox{if $*=3$}\\
		\R[\mbox{pt}] & \mbox{if $*=0$}\\
		0 & \mbox{otherwise}
	\end{array}\right. \]
Then the Mayer-Vietoris sequence involving the homology of $\cM^{\times 2}\setminus V^{\times 2}=(\cM\times(\cM\setminus V))\cup ((\cM\setminus V)\times \cM)$ is as follows.
\[ \begin{array}{c|cccccc}
	& & \scriptstyle(\cM\setminus V)^{\times 2} & & \scriptstyle\cM\times(\cM\setminus V) + (\cM\setminus V)\times\cM & & \scriptstyle\cM^{\times 2}\setminus V^{\times 2}\\\hline
	H_6 & \to & {\R[\partial\cM\otimes 1]+\R[1\otimes\partial\cM]}\atop{+\sum_{j,k}\R[\check a_j^i\otimes \check a_k^i]} 
		& \twoheadrightarrow & \R[\partial\cM\otimes 1]+\R[1\otimes\partial\cM] & \stackrel{0}{\to} & ?\\
	H_5 & \to & 0 & \to & 0 & \to & ?\\
	H_4 & \to & 0 & \to & 0 & \to & ?\\
	H_3 & \to & \sum_j(\R[1\otimes \check a_j^i]+\R[\check a_j^i\otimes 1]) & \hookrightarrow 
		& \sum_j(\R[1\otimes \check a_j^i]+\R[\check a_j^i\otimes 1]) & \stackrel{0}{\to} & ?\\
	H_2 & \to & 0 & \to & 0 & \to & ?\\
	H_1 & \to & 0 & \to & 0 & \to & ?\\
	H_0 & \to & \R & \to & \R+\R & \to & \R
	\end{array} \]
Here $\check a_j^i\eqdef a_j^i[4]$. Therefore the homology of $\cM^{\times 2}\setminus V^{\times 2}$ of dimensions at most 6 is
\[ H_*(\cM^{\times 2}\setminus V^{\times 2})=\left\{
	\begin{array}{ll}
		0 & \mbox{if $1\leq *\leq 6$}\\
		\R & \mbox{if $*=0$}
	\end{array}\right. \]
The homology of $C_2(M)\setminus C_2(V)$ is computed by the exact sequence:
\[ 	\to H_*(\breve{C}_2(\cM)\setminus \breve{C}_2(V)) 
	\to H_*(\cM^{\times 2}\setminus V^{\times 2}) 
	\to H_{*}(\cM^{\times 2}\setminus V^{\times 2}, \breve{C}_2(\cM)\setminus \breve{C}_2(V))
	\to \cdots\]
By excision, we have
\[ \begin{split}
	H_{*}(\cM^{\times 2}\setminus V^{\times 2}, \breve{C}_2(\cM)\setminus \breve{C}_2(V))
	&\cong
	H_{*}((\cM\setminus V)\times\R^7, (\cM\setminus V)\times(\R^7\setminus \{0\}))\\
	&\cong
	H_{*-7}(\cM\setminus V)\otimes H_6(S^6).
\end{split} \]
In particular, $H_{*}(\cM^{\times 2}\setminus V^{\times 2}, \breve{C}_2(\cM)\setminus \breve{C}_2(V))=0$ for $0\leq *\leq 6$. Thus the above exact sequence is as follows.
\[ \begin{array}{c|cccccc}
	& & \scriptstyle\breve{C}_2(\cM)\setminus\breve{C}_2(V)\ & & \scriptstyle\cM^{\times 2}\setminus V^{\times 2} & & \scriptstyle(\cM^{\times 2}\setminus V^{\times 2}, \breve{C}_2(\cM)\setminus \breve{C}_2(V))\\\hline
	H_6 & \to & ? & \to & 0 & \to & 0\\
	H_5 & \to & 0 & \to & 0 & \to & 0\\
	H_4 & \to & 0 & \to & 0 & \to & 0\\
	H_3 & \to & 0 & \to & 0 & \to & 0\\
	H_2 & \to & 0 & \to & 0 & \to & 0\\
	H_1 & \to & 0 & \to & 0 & \to & 0\\
	H_0 & \to & \R & \to & \R & \to & 0
	\end{array} \]

Then the homology of $C_2(M)\setminus C_2(V)\cup STV$ is computed as follows. Note that this space can be obtained by gluing $ST\cM\cong \cM\times S^6$ and $C_2(M)\setminus C_2(V)$ along $ST(\cM\setminus V)\cong (\cM\setminus V)\times S^6$. The Mayer-Vietoris sequence is as follows.
\[ \begin{array}{c|cccccc}
	& & \scriptstyle(\cM\setminus V)\times S^6 & & \scriptstyle\cM\times S^6+C_2(M)\setminus C_2(V) & & \scriptstyle C_2(M)\setminus C_2(V)\cup STV\\\hline
	H_5 & \to & 0 & \to & 0 & \to & 0\\
	H_4 & \to & 0 & \to & 0 & \to & ?\\
	H_3 & \to & \sum_j\R[\check a_j^i\otimes 1] & \to & 0 & \to & 0\\
	H_2 & \to & 0 & \to & 0 & \to & 0\\
	H_1 & \to & 0 & \to & 0 & \to & 0\\
	H_0 & \to & \R & \to & \R+\R & \to & \R
	\end{array} \]
Hence $H_*(D_{\{i\}}(\beta_M^0))$ vanishes at $*=1,2,3,5$. This shows that $E^2_{p,q}= 0$ if $p+q\leq 6$ and ($p\notin\{0,2\}$ or $q\notin\{0,4,6\}$). Moreover, all differentials $E^2_{*,*}\to E^2_{*-2,*+1}$ involving $E^2_{p,q}$ ($p+q\leq 6$) are zero and hence $E^2_{p,q}=E^\infty_{p,q}$ there.
\end{proof}

\subsection*{Lescop cycles $F(a)$}

In order to prove Lemma~\ref{lem:E24}, we shall give a higher dimensional analogue of the Lescop cycles, which were constructed by Lescop in 3-dimension \cite{Les2}. Namely, for each $a=a_j^i$, we consider a 6-cycle $F(a)$ on the configuration space bundle $D(E^\Gamma(i))$ of the form:
\[ \begin{split}
	F(a)\eqdef &(C(a)\times S^2)\\
	&\cup -(\widetilde{S}(a)\widetilde\times (4\times p(a)))
	\cup -((4\times p(a))\widetilde\times \widetilde{S}(a))\\
	&\cup \mathrm{diag}(n)(\widetilde{S}(a)) \\
	& \mbox{($p(a)$: base point of $a$)}
	\end{split} \]
where $\widetilde{S}(a)\widetilde\times (4\times p(a))\eqdef\cup_t\{x_t\times (4\times p(a)_t)\,|\, x_t\in \widetilde{S}(a)_t\}$ and $(4\times p(a))\widetilde\times\widetilde{S}(a)$ is its symmetric. In order to define $F(a)$, we choose a vector field $n$ that is a section of the trivial $S^6$-bundle $\widetilde{STV}[4]$ (the sub $STV[4]$-bundle of $\tilde{V}[4]$) restricted to $\widetilde{S}(a)$ such that near $\partial\widetilde{S}(a)$ it is normal to $\widetilde{S}(a)$ and tangent to $\partial \tilde{V}[4]$. Moreover we assume that the map
\begin{equation}\label{eq:mapping-deg}
 (\widetilde{S}(a), \partial\widetilde{S}(a))\to (S^6,*) 
\end{equation}
given by the trivialization composed with the projection to the $S^6$-factor and by $n$, is mapping degree 0 so that $F(a)$ represents a class in $E_{2,4}^\infty$. Then we introduce a local coordinate $a\times [0,1]\subset \partial V$ where the second coordinate determined by the direction of $n$.

\begin{figure}
\psfrag{T(a,a)}[cc][cc]{$\scriptstyle T(0\times a\times 0,0\times a\times 1)$}
\psfrag{A(0,1)}[cc][cc]{$\scriptstyle A(0,1)$}
\psfrag{0a0}[cc][cc]{$\tiny{a\times[-([0,4]\times p(a)\times 1)}\atop{\cup(4\times p(a)\times [0,1])]}$}
\psfrag{4p0}[cc][cc]{$\tiny {(4\times p(a)\times 0)}\atop{(0\times a\times [0,1])}$}
\psfrag{[04]p0}[cc][cc]{$\tiny {([0,4]\times p(a)\times 0)}\atop{(0\times a\times 1)}$}
\psfrag{S4p}[cc][cc]{$\scriptstyle -\widetilde{S}(a)\widetilde{\times}(4\times p(a));$}
\psfrag{diag}[cc][cc]{$\scriptstyle \mathrm{diag}(n)(\widetilde{S}(a))$}
\psfrag{d}[cc][cc]{$\partial$}
\fig{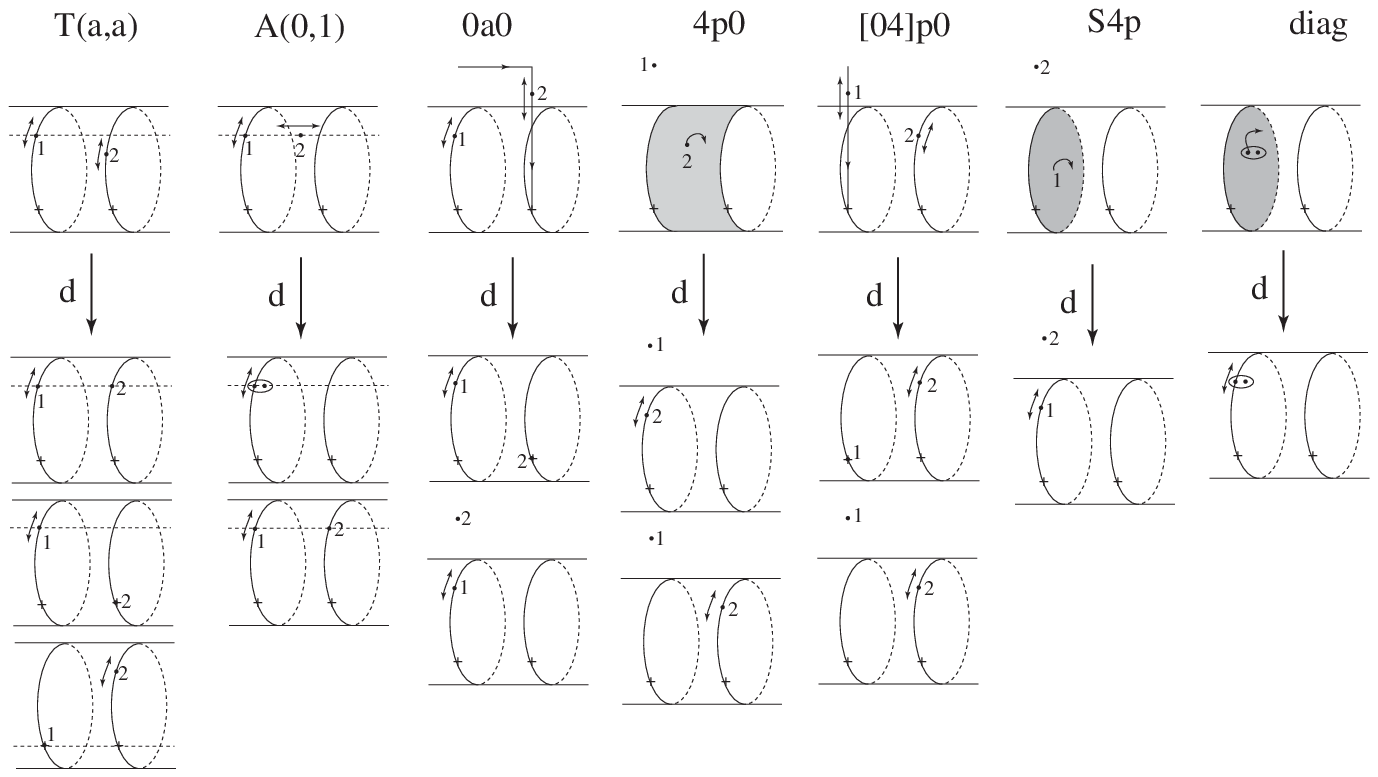}
\caption{Lescop cycle $F(a)$}\label{fig:F(a)}
\end{figure}

The $C(a)$ is a 4-chain on $C_2([0,4]\times a\times [0,1])\subset C_2(M_t)\setminus C_2(V[-1]_t)\cup STV_t$ defined as a sum of the following chains:
\begin{itemize}
\item $T(0\times a\times 0,0\times a\times 1)$
\item $A(0,1)$
\item $(0\times a\times 0)\times\bigl[-([0,4]\times p(a)\times 1)\cup(4\times p(a)\times[0,1])\bigr]$
\item $\bigl((4\times p(a)\times 0)\times (0\times a\times [0,1])\bigr)
	\cup \bigl(([0,4]\times p(a)\times 0)\times (0\times a\times 1)\bigr)$
\end{itemize}
To describe $T(0\times a\times 0,0\times a\times 1)$ and $A(0,1)$, we identify $S^3$ with $\R^3\cup\{\infty\}$ and we consider the 4-dimensional submanifold $T$ of $B\ell(S^3,\{\infty\})\times S^3$ defined by
\[ 	T\eqdef\overline{
	\{(x,y,z)\times(x,y,z')\,|\,
		x,y,z,z'\in\R,\ 
		z\geq z'\}
	}
	\subset B\ell(S^3,\{\infty\})\times S^3 \]
with
\[ \partial T=\overline{\{(x,y,z)\times (x,y,z)\}}
	\cup\overline{\{\{\infty\}\times (x,y,z)\}}
	\cup\overline{\{(x,y,z)\times\{\infty\}\}}. \]
Consider a pair of parallel cycles $0\times a\times 0$ and $0\times a\times 1$ and identify $(0\times a\times 0)\times (0\times a\times 1)$ by the base point preserving ($p(a)\leftrightarrow \{\infty\}$) diffeomorphism
\[ \varphi: (0\times a\times 0)\times (0\times a\times 1)\to S^3\times S^3\subset B\ell(S^3,\{\infty\})\times S^3. \]
Then we set
\[ \begin{split}
	T(0\times a\times 0,0\times a\times 1)&\eqdef \varphi^{-1}T\\
	A(0,1)&\eqdef \{(x\times 0)\times (x\times s)\,|\,x\in a, s\in[0,1)\}\\
	&\subset (a\times 0)\times (a\times [0,1])
	\end{split} \]

The chain $\mathrm{diag}(n)(\widetilde{S}(a))$ denotes the image of $\widetilde{S}(a)$ in the trivial $S^6$-bundle $\widetilde{STV}[4]$ under the section $n$. See Figure~\ref{fig:F(a)} for the form of $F(a)$. Lemma~\ref{lem:E24} follows from Lemma~\ref{lem:F(a)-null} and \ref{lem:beta-vanish} described in the following.

\begin{Lem}\label{lem:F(a)-null}
\begin{enumerate}
\item $[F(a)]$ spans $E_{2,4}^\infty(D(E^\Gamma(i)))$.
\item $F(a)$ is null in $H_6(C(E^\Gamma(i));\R)$.
\end{enumerate}
\end{Lem}
\begin{proof}
According to the proof of Lemma~\ref{lem:ssd} and from the definition of $F(a)$, the image of $[F(a)_t]$ under the Mayer-Vietoris boundary homomorphism is $[\check{a}\otimes 1]$ and moreover its collection over the $S^2$ is $[\check{a}\otimes S^2]$ in $H_6(D(E^\Gamma(i)))$. Hence $[F(a)]$ spans $E_{2,4}^\infty(D(E^\Gamma(i)))$.

The second assertion follows from the naturality of the Serre spactral sequences (see e.g., \cite{Hat}), in our setting together with Lemma~\ref{lem:ssd} implying that there are homomorphisms between $E_{*,*}^\infty$'s induced by the inclusion
\[ \begin{split}
	E_{0,6}^\infty(D(E^\Gamma(i))) & \to E_{0,6}^\infty(C(E^\Gamma(i)))\\
	E_{2,4}^\infty(D(E^\Gamma(i))) & \to E_{2,4}^\infty(C(E^\Gamma(i)))=0
	\end{split} \]
which is isomorphism on $E_{0,6}^\infty$ and is zero map on $E_{2,4}^\infty$.
\end{proof}

\begin{Lem}\label{lem:beta-vanish}
The 6-form $\beta_{M_t}^0$ on $D(E^\Gamma(i))$ evaluated on any cycle of $E_{2,4}^\infty(D(E^\Gamma(i)))$ vanishes.
\end{Lem}
\begin{proof}
We prove that
\[ \int_{F(a)}\beta_{M_t}^0=0. \]

First extend the form $\beta_M$ on $C_2(M)$ obviously to a fundamental 6-form on the trivial bundle $C_2(M)\times S^2$ and denote it also by $\beta_M$. We have $\int_{C(a)\times S^2}\beta_{M_t}^0=\int_{C(a)\times S^2}\beta_{M}=0$ since $C(a)$ lives inside $C_2([0,4]\times a\times [0,1])\subset C_2(M_t)$ where $\beta_{M_t}^0$ and $\beta_M$ coincide.

The normalization of Proposition~\ref{prop:normalize-part} implies that the integrals vanish on
\[ -(\widetilde{S}(a)\widetilde\times (4\times p(a)))\cup
	-((4\times p(a))\widetilde\times \widetilde{S}(a)). \]

Since $F(a)$ is null homologous in $C(E^\Gamma(i))$ by Lemma~\ref{lem:F(a)-null}, it is enough to prove that
\[ \int_{\mathrm{diag}(n_0)(\widetilde{S}_0(a))}\beta_M=\int_{\mathrm{diag}(n)(\widetilde{S}(a))}\beta_{M_t}^0 \]
where $\widetilde{S}_0(a)$ is any embedding of a 6-manifold diffeomorphic to $\widetilde{S}(a)$ into the trivial sub bundle $V[4]\times S^2$ of $M\times S^2$ having the same behavior as $\widetilde{S}(a)$ near $\partial V\times S^2$, and $n_0$ is any vector field on $\widetilde{S}_0(a)$ tangent to the fibers of $V[4]\times S^2$ which coincides with $n$ near the boundary and which satisfies the same constraint as $n$ on mapping degree of the map (\ref{eq:mapping-deg}). Then the boundary relative homology classes of the images of the sections $n$ and $n_0$ coincide and hence the integrals also coincide.
\end{proof}

\end{appendix}

%%%%%%%%%%%%%%%%%%%%%%%%%%%%%%
%%%%%%%%%%%%%%%%%%%%%%%%%%%%%%
%%%%%%%%%%%%%%%%%%%%%%%%%%%%%%
%%%%%%%%%%%%%%%%%%%%%%%%%%%%%%%%%%%%%%%

\end{document}